\documentclass[12pt]{article}
\usepackage{graphicx}
\usepackage{amsmath,amsthm,amssymb,enumerate}%, esint}
\usepackage{euscript,mathrsfs}
\usepackage{color}
\usepackage{dsfont}
\usepackage[left=2cm,right=2cm,top=3.5cm,bottom=3.5cm]{geometry}
\usepackage{color}
\usepackage{subfigure}
\usepackage[framemethod=tikz]{mdframed}
\allowdisplaybreaks

\usepackage{soul}
\usepackage{bm}

\catcode`\@=11 \@addtoreset{equation}{section}

\catcode`\@=12

\newtheorem{Theorem}{Theorem}[section]
\newtheorem{Proposition}[Theorem]{Proposition}
\newtheorem{Lemma}[Theorem]{Lemma}
\newtheorem{Corollary}[Theorem]{Corollary}

\theoremstyle{definition}
\newtheorem{Definition}[Theorem]{Definition}

\newtheorem{Remark}[Theorem]{Remark}

\newcommand{\bTheorem}[1]{
	\begin{Theorem} \label{T#1} }
	\newcommand{\eT}{\end{Theorem}}

\newcommand{\bProposition}[1]{
	\begin{Proposition} \label{P#1}}
	\newcommand{\eP}{\end{Proposition}}

\newcommand{\bLemma}[1]{
	\begin{Lemma} \label{L#1} }
	\newcommand{\eL}{\end{Lemma}}

\newcommand{\bCorollary}[1]{
	\begin{Corollary} \label{C#1} }
	\newcommand{\eC}{\end{Corollary}}

\newcommand{\bRemark}[1]{
	\begin{Remark} \label{R#1} }
	\newcommand{\eR}{\end{Remark}}

\newcommand{\bDefinition}[1]{
	\begin{Definition} \label{D#1} }
	\newcommand{\eD}{\end{Definition}}

\newcommand{\Td}{\mathbb{T}^d}

\newcommand{\prst}{\mathbb{P}}

\newcommand{\bfomega}{\boldsymbol{\omega}}

\newcommand{\jump}[1]{\left[ \left[ #1 \right] \right]}

\newcommand{\vrh}{\vr_h}

\newcommand{\tvm}{\widetilde{\vc{m}}}

\newcommand{\bfphi}{\boldsymbol{\varphi}}

\newcommand{\bFormula}[1]{
	\begin{equation} \label{#1}}
	\newcommand{\eF}{\end{equation}}

\newcommand{\grid}{\mathcal{T}}
\newcommand{\facesK}{\mathcal{E}(K)}
\newcommand{\vuh}{\vu_h}

\newcommand{\TS}{\Delta t}

\newcommand{\Divh}{{\rm div}_h}

\newcommand{\Ov}[1]{\overline{#1}}

\newcommand{\vr}{\varrho}

\newcommand{\tvr}{\tilde \vr}

\newcommand{\vu}{\vc{u}}
\newcommand{\vm}{\vc{m}}

\newcommand{\vn}{\vc{n}}

\newcommand{\vc}[1]{{\bm #1}}

\newcommand{\Div}{{\rm div}_x}
\newcommand{\Grad}{\nabla_x}

\newcommand{\dx}{\,{\rm d} {x}}

\newcommand{\dt}{\,{\rm d} t }

\newcommand{\vU}{\vc{U}}

\newcommand{\intG}[1]{\int_{\sigma} #1 \ {\rm dS}_x}

\newcommand{\sumG}{\sum_{\sigma \in \mathcal{E}}}

\newcommand{\intTd}[1]{\int_{\mathbb{T}^d} #1 \ \dx}

\newcommand{\D}{{\rm d}}

\newcommand{\expe}[1]{ \mathbb{E} \left[ #1 \right] }

\newcommand{\br}{ \nonumber \\ }

\def\softd{{\leavevmode\setbox1=\hbox{d}%
		\hbox to 1.05\wd1{d\kern-0.4ex{\char039}\hss}}}
\newcommand{\Fup}{{F_h^\up}}
\newcommand{\vFup}{{\vc{F}_h^\up}}
\newcommand{\up}{ {\rm up} }
\newcommand{\avs}[1]{\left\{\left\{ #1\right\}\right\}}
\renewcommand{\avs}[1]{\left\{\!\!\left\{ #1\right\}\!\!\right\}}
\newcommand{\Gradedge}{\nabla_{\mathcal{D}}}
\newcommand{\bfpphi}{\boldsymbol{\phi}}

\definecolor{Cgrey}{rgb}{0.85,0.85,0.85}
\definecolor{Cblue}{rgb}{0.50,0.85,0.85}
\definecolor{Cred}{rgb}{1,0,0}
\definecolor{fancy}{rgb}{0.10,0.85,0.10}

\newcommand\Cbox[2]{%
	\newbox\contentbox%
	\newbox\bkgdbox%
	\setbox\contentbox\hbox to \hsize{%
		\vtop{
			\kern\columnsep
			\hbox to \hsize{%
				\kern\columnsep%
				\advance\hsize by -2\columnsep%
				\setlength{\textwidth}{\hsize}%
				\vbox{
					\parskip=\baselineskip
					\parindent=0bp
					#2
				}%
				\kern\columnsep%
			}%
			\kern\columnsep%
		}%
	}%
	\setbox\bkgdbox\vbox{
		\color{#1}
		\hrule width  \wd\contentbox %
		height \ht\contentbox %
		depth  \dp\contentbox
		\color{black}
	}%
	\wd\bkgdbox=0bp%
	\vbox{\hbox to \hsize{\box\bkgdbox\box\contentbox}}%
	\vskip\baselineskip%
}

\mdfdefinestyle{MyFrame}{%
	linecolor=black,
	outerlinewidth=1pt,
	roundcorner=5pt,
	innertopmargin=\baselineskip,
	innerbottommargin=\baselineskip,
	innerrightmargin=10pt,
	innerleftmargin=10pt,
	backgroundcolor=white!20!white}

\begin{document}

%%%%%%%%%%%%%%%%%%%%%%%%%%%%%%%%

\title{Monte Carlo method and the random \\ isentropic Euler system}

\author{Eduard Feireisl
	\thanks{The work of E.F. and H.M. was partially supported by the
		Czech Sciences Foundation (GA\v CR), Grant Agreement
		24-11034S. The Institute of Mathematics of the Academy of Sciences of
		the Czech Republic is supported by RVO:67985840. \newline
		\hspace*{1em} $^\spadesuit$
		The work of M.L. and  C.Y. was funded by the Deutsche Forschungsgemeinschaft (DFG, German Research Foundation) - Project number 525853336 -- SPP 2410 ``Hyperbolic Balance Laws: Complexity, Scales and Randomness.''  M.L. is grateful to the Gutenberg Research College
		and Mainz Institute of Multiscale Modelling for supporting her research.
	} \and  M\'aria Luk\'a\v{c}ov\'a-Medvi{\softd}ov\'a$^{\spadesuit}$
	\and  Hana Mizerov\'a $^{ *,\dagger}$
\and Changsheng Yu $^{\spadesuit}$
}

\date{\today}

\maketitle

\bigskip

\centerline{$^*$  Institute of Mathematics of the Academy of Sciences of the Czech Republic}

\centerline{\v Zitn\' a 25, CZ-115 67 Praha 1, Czech Republic}

\medskip

\centerline{$^\spadesuit$ Institute of Mathematics, Johannes Gutenberg--University Mainz}

\centerline{Staudingerweg 9, 55 128 Mainz, Germany}

\medskip

\centerline{$^\dagger$ Department of Mathematical Analysis and Numerical Mathematics, Comenius University}
\centerline{Mlynsk\' a dolina, 842 48 Bratislava, Slovakia}

\begin{abstract}
	
We show several results on convergence of the Monte Carlo method applied
to consistent approximations of the isentropic Euler system of gas dynamics with uncertain initial data.
Our method is based on combination of several new concepts. We work with the dissipative weak solutions that can be seen as a universal closure of consistent approximations. Further, we apply the set-valued version of the Strong law of large numbers for general multivalued mapping with closed range and the Koml\'os theorem on strong converge of empirical averages of integrable functions. Theoretical results are illustrated by a series of numerical simulations obtained by an unconditionally convergent viscosity finite volume method combined with the Monte Carlo method.
%\begin{itemize}
%	\item
%	
%	The concept of dissipative weak solutions as a universal closure of families of consistent approximations.
%	
%	\item Set--valued version of the Strong law of large numbers for general multivalued mapping with closed
%	range.
%	
%	\item Application of Koml\' os version of convergence of empirical averages of integrable functions.
%	
%	\item Unconditionally convergent finite volume approximation scheme.
%	
%\end{itemize} 	
%
%	
%	
\end{abstract}

%\bigskip


%\tableofcontents

\section{Introduction}
\label{E}

The Euler system of gas dynamics is a standard example of a hyperbolic system of nonlinear conservation
laws. We consider its simplified isentropic version describing the time evolution of the
mass density $\vr = \vr(t,x)$ and the velocity $\vu = \vu(t,x)$ of a compressible fluid:
\begin{align}
	\partial_t \vr + \Div (\vr \vu) &= 0, \label{E1} \\
	\partial_t (\vr \vu) + \Div (\vr \vu \otimes \vu) + a \Grad \vr^\gamma  &= 0,\ a > 0,\ \gamma > 1
	\label{E2}
\end{align}
The model ignores the effect of temperature changes, or, more precisely,
imposes the isentropic regime, where the pressure is directly related to the density.
For the sake of simplicity, we focus on the space periodic boundary conditions,
\begin{equation} \label{E3}
	x \in \Td = \left( [-\pi, \pi]|_{\{ - \pi, \pi \} } \right)^d,\ d  = 2,3.
	\end{equation}
The problem is formally closed by prescribing the initial state
\begin{equation} \label{E4}
	\vr(0, \cdot) = \vr_0,\ \vr \vu(0, \cdot) = \vr_0 \vu_0 \equiv \vm_0.
	\end{equation}

\subsection{Well posedness in the class of smooth solutions}

As is well known, the isentropic Euler system can be written as a symmetric hyperbolic system in the
variables $r = \vr^{\frac{\gamma - 1}{2}}$, $\vu$. Accordingly, the abstract theory developed in the
monograph of Majda \cite{Majd} applies yielding the following result.

\begin{Proposition} [{\bf Local existence of smooth solutions}] \label{PE1}
Suppose
\begin{equation} \label{E5}
\vr_0 \in W^{k,2}(\Td),\ \inf_{x \in \Td} \vr_0(x) > 0,\ \vu_0 \in W^{k,2}(\Td; R^d),\ k > 1 + \frac{d}{2}.
\end{equation}

Then there exists $0 < T_{\rm max} \leq \infty$	and a classical solution $(\vr, \vu)$ of the Euler system \eqref{E1}--\eqref{E4}
defined in $[0, T_{\rm max})$, unique in the class
\begin{equation} \label{E6}
\vr \in C_{\rm loc}([0, T_{\rm max}); W^{k,2}(\Td)), \ \vu \in C_{\rm loc}([0, T_{\rm max}); W^{k,2}(\Td; R^d)).
\end{equation}
In addition,
\begin{equation} \label{E7}
	T_{\rm max} < \infty \ \Rightarrow \ \lim_{t \to T_{\rm max} - } \Big( \| \vr(t, \cdot) \|_{W^{k,2}(\Td)} +
	\| \vu(t, \cdot) \|_{W^{k,2}(\Td;R^d)} +  \left\| \frac{1}{\vr (t, \cdot) } \right\|_{C(\Td)} \Big) = \infty.
\end{equation}

\end{Proposition}

\subsection{Weak solutions}

Unfortunately, the smooth solutions may experience singularities in the form of shock waves that develop in a finite time. Extending
smooth solutions beyond the blow--up time $T_{\rm max}$ requires a weaker concept of solutions, where derivatives are interpreted as
distributions. The weak solutions, however, are not uniquely
determined by the initial data, for standard examples see, e.g., the
monograph by Smoller \cite{SMO}.
To save, at least formally, the desired well posedness of the problem, the weak formulation of the field
equations \eqref{E1}, \eqref{E2} is usually accompanied by a variant of energy balance (``entropy'' inequality). To this end,
it is more convenient to rewrite the system in terms of the conservative variables $\vr$, and the momentum $\vm = \vr \vu$. The
associated total energy reads
\[
E(\vr, \vm) = \frac{1}{2} \frac{|\vm|^2}{\vr} + \frac{a}{\gamma - 1} \vr^\gamma.
\]
More precisely, we define
\begin{equation} \label{E13}
	E(\vr, \vm) = \left\{ \begin{array}{l} \frac{1}{2} \frac{|\vm|^2}{\vr} + \frac{a}{\gamma - 1} \vr^\gamma \ \mbox{whenever}\ \vr > 0,  \\  0 \ \mbox{if}\ \vr = 0,\ \vm = 0, \\ \infty \ \mbox{otherwise}                        \end{array} \right.
\end{equation}
Accordingly,  $E: R^{d + 1} \to [0, \infty]$ is a convex l.s.c function, strictly convex in the interior of its domain.
\begin{Definition}[{\bf Admissible weak solution}] \label{DE1}
	Let the initial data satisfy
	\begin{equation} \label{E14}
		\vr_0 \in L^1(\Td), \ \vm_0 \in L^1(\Td; R^d),\ \intTd{ E(\vr_0, \vm_0) } < \infty.
	\end{equation}
	
	We say that $(\vr, \vm)$ is an \emph{admissible weak solution} to the Euler system in $[0,T) \times \Td$, $0 < T \leq \infty$, if the following holds:
	\begin{itemize}
		\item {\bf Regularity.}
		The solution $(\vr, \vm)$ belongs to the class
		\begin{align}
			\vr \in C_{\rm weak, loc}([0,T); L^\gamma(\Td)),\
			\vm &\in C_{\rm weak, loc}([0,T); L^{\frac{2\gamma}{\gamma + 1}}(\Td; R^d)), \br
			\intTd{ E(\vr, \vm)(t, \cdot) } &\leq \intTd{ E(\vr_0, \vm_0) } \ \mbox{for any}\ 0 \leq t < T.
			\label{E15}
		\end{align}
		\item {\bf Equation of continuity.}
		The integral identity
		\begin{equation} \label{E16}
			\int_0^T \intTd{ \Big[ \vr \partial_t \varphi + \vm \cdot \Grad \varphi \Big] } \dt = - \intTd{ \vr_0 \varphi(0, \cdot)}
		\end{equation}
		holds for any $\varphi \in C^1_c([0,T) \times \Td)$.
		
		\item {\bf Momentum equation.} The integral identity
		\begin{equation} \label{E17}
			\int_0^T \intTd{ \left[ \vm \cdot \partial_t \bfphi + \mathds{1}_{\vr > 0} \frac{\vm \otimes \vm}{\vr} : \Grad \bfphi +
				a \vr^\gamma \Div \bfphi \right] } \dt = - \intTd{ \vm_0 \cdot \bfphi(0, \cdot) }
		\end{equation}
		holds for any $\bfphi \in C^1_c([0,T) \times \Td; R^d)$.
		
		\item {\bf Energy inequality.}
		The integral inequality
		\begin{align}
			\int_0^T \intTd{ \left[ E(\vr, \vm) \partial_t \varphi + \mathds{1}_{\vr > 0} \Big( E(\vr, \vm) + a \vr^\gamma \Big) \frac{\vm}{\vr} 	
				\cdot \Grad \varphi \right] } \dt \geq - \intTd{ E(\vr_0, \vm_0) \varphi(0, \cdot) }
			\label{E18}
		\end{align}
		holds for any $\varphi \in C^1_c([0,T) \times \Td)$, $\varphi \geq 0$.
		
	\end{itemize}
	
\end{Definition}

\begin{Remark} \label{reg}
	The reader will have noticed that \eqref{E18} already includes
	the inequality in \eqref{E15}.
	\end{Remark}

The recent applications of the theory of convex integration to problems in fluid dynamics revealed
a number of rather disturbing facts concerning well posedness of certain problems, including the isentropic Euler
system in higher space dimensions $d=2,3$, see Buckmaster et al. \cite{BuDeSzVi} \cite{BucVic},
Chiodaroli et al. \cite{Chiod}, \cite{ChiDelKre}, De Lellis and Sz\' ekelyhidi \cite{DelSze3}, or, more recently Giri and Kwon
\cite{GirKwo}, to name only a few. The name \emph{wild data} is used for the initial state that gives rise to infinitely
many \emph{admissible} weak solutions on an arbitrarily short time interval.
\begin{Definition} [{\bf Wild data}] \label{DE2}
	
	We say that the initial data $(\vr_0, \vm_0)$ are \emph{wild}, if there exists $T > 0$ such that the Euler system admits infinitely many
	admissible weak solutions in $[0,T] \times \Td$ distinct on any time interval $[0,\tau]$, $0 < \tau \leq T$.

\end{Definition}

The following result was proved in \cite[Theorem 1.3]{ChFe2022}.

\begin{Proposition} [{\bf Density of wild data}] \label{PE3}
	Let $d=2,3$.
	The set of wild data $(\vr_0, \vm_0)$ is dense in the space $L^\gamma (\Td; [0, \infty)) \times L^{\frac{2 \gamma}{\gamma + 1}}
	(\Td; R^d)$.	
\end{Proposition}

\begin{Remark} \label{RE2}
	The wild data give rise to a family of admissible solutions that are, in general, only local in time. If we drop the
	differential version of the energy inequality \eqref{E18} and retain only \eqref{E15}, the density of the corresponding wild solutions that give rise to \emph{global in time} weak solutions
	was shown by Chen, Vasseur, and Yu \cite{ChVaYu}.
	
\end{Remark}

\begin{Remark} \label{RE3}
	
	Adapting the method of Glimm, quite regular wild data can be constructed, in fact piecewise smooth out of a finite number of straight--lines,
	with the associated weak solutions with the density bounded below away from zero, see
	\cite{ChFe2023}.
	
\end{Remark}

\subsection{Dissipative weak solutions}

Despite the ``positive'' existence result stated in Proposition \ref{PE3}, the existence of admissible weak solutions for \emph{arbitrary} finite energy initial data is not known. Moreover, in view of the ill--posedness discussed in the preceding part,
the relevance of the Euler system to describe the behaviour of fluids in higher
space dimensions may be dubious. Indeed the Euler system is a model of a perfect (inviscid) fluid and as such
should reflex the behaviour of real (viscous) fluids in the vanishing viscosity limit. The low viscosity
regime, however, is characteristic for turbulence, where the solutions may develop fast oscillations. Mathematically speaking, they may approach the limit state only in the weak sense. Unfortunately,
in view of the arguments developed in \cite{MarEd}, weak limits of the compressible Navier--Stokes system
(on the whole physical space $R^d$) cannot be weak solutions of the Euler system.

The concept of \emph{dissipative weak solution}, developed in the numerical context in \cite{FeLMMiSh}, reflects
better the idea of understanding the Euler system as a limit of \emph{consistent approximations}, among which
the vanishing viscosity limit.

Let $\mathcal{M}^+(\Td, R^{d \times d}_{\rm sym})$ denote the set of all non--negative matrix valued measures on $\Td$.

\begin{Definition}[{\bf Dissipative weak solution}] \label{DE3}
	Let the initial data satisfy
	\[
	\vr_0 \in L^1(\Td), \ \vm_0 \in L^1(\Td; R^d),\ \intTd{ E(\vr_0, \vm_0) } < \infty.
	\]
	We say that $(\vr, \vm)$ is \emph{dissipative weak (DW) solution} to the Euler system in $[0,T) \times \Td$,
	$0 < T \leq \infty$, if the following holds:
	\begin{itemize}
		
		\item {\bf Regularity.}
		The solution $(\vr, \vm)$ belongs to the class
		\begin{align}
			\vr &\in C_{\rm weak, loc}([0,T); L^\gamma(\Td)),\
			\vm \in C_{\rm weak, loc}([0,T); L^{\frac{2\gamma}{\gamma + 1}}(\Td; R^d)),\br
			\intTd{ E(\vr, \vm)(t, \cdot) } &\leq \intTd{ E(\vr_0, \vm_0) } \ \mbox{for any}\ 0 \leq t < T.
			\nonumber
		\end{align}
		\item {\bf Equation of continuity.}
		The integral identity
		\[
		\int_0^T \intTd{ \Big[ \vr \partial_t \varphi + \vm \cdot \Grad \varphi \Big] } \dt =  - \intTd{ \vr_0 \varphi(0, \cdot)}
		\]
		holds for any $\varphi \in C^1_c([0,T) \times \Td)$.
		
		\item {\bf Momentum equation.}  The integral identity
		\begin{align}
			\int_0^T &\intTd{ \left[ \vm \cdot \partial_t \bfphi + \mathds{1}_{\vr > 0} \frac{\vm \otimes \vm}{\vr} : \Grad \bfphi +
				a \vr^\gamma \Div \bfphi \right] } \dt \br &=
			- \int_0^T \intTd{ \Grad \bfphi : \D \mathfrak{R} (t) }  - \intTd{ \vm_0 \cdot \bfphi(0, \cdot) }
			\label{E19}
		\end{align}
		for any $\bfphi \in C^1_c([0,T) \times \Td; R^d)$, where
		\begin{equation} \label{E20}
			\mathfrak{R} \in L^\infty([0, T); \mathcal{M}^+(\Td; R^{d \times d}_{\rm sym})).	
		\end{equation}
		
		\item{\bf Compatibility, energy defect.} There exists a non--increasing function $\mathcal{E}:
		[0,T) \to [0, \infty)$ satisfying
		\begin{align}
			\mathcal{E}(0-) &= \intTd{ E(\vr_0, \vm_0) }, \br
			\mathcal{E}(\tau+) - \intTd{ E(\vr, \vm)(\tau, \cdot) } &\geq \min \left\{ \frac{1}{2}, \frac{1}{\gamma - 1} \right\} \int_{\Td} \D ({\rm trace} [\mathfrak{R}] (\tau)) 	
			\label{E21}	
		\end{align}
		for any $0 \leq \tau < T$.	
	\end{itemize}
	
\end{Definition}

Obviously, there is certain freedom in the choice of the function $\mathcal{E}$. To pick up a physically
relevant solution, we consider the class of maximal (DW) solutions.

\begin{Definition}[{\bf Maximal (DW) solution}] \label{DE4}
	
	A (DW) solution $(\vr, \vm)$ with the associated function $\mathcal{E}$ on a time interval $[0,T)$ is
	called \emph{maximal} if the following implication holds:
	
	Suppose there is another (DW) solution $(\tvr, \tvm)$ of the same problem with the associated function $\widetilde{\mathcal{E}}$
	such that
	\[
	\widetilde{\mathcal{E}} \leq \mathcal{E} \ \mbox{in}\ [0,T).
	\]
	Then
	\[
	\widetilde{\mathcal{E}} = \mathcal{E}.
	\]
	
\end{Definition}

Unlike the admissible weak solutions, the (DW) solutions are known to exist globally in time for any finite energy initial data. The following result was proved in \cite[Theorem 2.5]{BreFeiHof19}.

\begin{Proposition}[{\bf Global existence of maximal (DW) solution}] \label{PE4}
	For any measurable initial data $(\vr_0, \vm_0)$,
	\[
	 \intTd{ E(\vr_0, \vm_0) } < \infty,
	\]
	the Euler system admits a maximal (DW) solution $(\vr, \vm)$ in $[0, \infty) \times \Td$. Moreover, the
	solution can be selected in such a way that the mapping
	\[
	(\vr_0, \vm_0) \in W^{-\ell,2}(\Td) \times W^{-\ell,2}(\Td; R^d)
	\mapsto (\vr, \vm) \in C_{{\rm loc}}([0, \infty); W^{-\ell,2}(\Td) \times
	W^{-\ell, 2}(\Td; R^d))
	\]
	is Borel measurable for any $\ell > d$.
	
\end{Proposition}

The energy defect of the maximal (DW) solutions approaches zero as $t \to \infty$, see \cite[Theorem 2.3]{Fei2020}.

\begin{Proposition}[{\bf Vanishing energy defect}] \label{PE5}
	Let $(\vr, \vm)$ be a maximal (DW) solution of the Euler system in $[0,\infty) \times \Td$.
	
	Then the associated energy function $\mathcal{E}$ satisfies
	\[
	\lim_{t \to \infty} \mathcal{E}(t) = \lim_{t \to \infty} \intTd{ E(\vr, \vm) (t,\cdot) },
	\]
	in particular
	\[
	{\rm ess} \sup_{t > \tau} \| \mathfrak{R}(t) \|_{\mathcal{M}^+(\Td; R^{d \times d}_{\rm sym})} \to 0
	\ \mbox{as} \ \tau \to \infty .
	\]	
\end{Proposition}

Finally, we report the following weak--strong uniqueness result, see \cite[Theorem 5.2]{FeiGhoJan}.

\begin{Proposition}[{\bf Weak--strong uniqueness}] \label{PE6}
	Let $(\vr, \vm)$ be a (DW) solution of the Euler system in $[0,T) \times \Td$, $0 < T < \infty$ in the sense of Definition \ref{DE3}. Suppose that
	$(\tvr, \tvm)$ is a weak solution of the same problem in $[0,T) \times \Td$ satisfying:
	\begin{itemize}
		\item $\tvm = \tvr \vU$, where
		\[
		\tvr \in B^{\alpha, \infty}_p ([\delta,T] \times \Td) \cap C([0,T]; L^1(\Td)), \
		\vU \in B^{\alpha, \infty}_p ([\delta,T] \times \Td; R^d) \cap C([0,T]; L^1(\Td; R^d))
		\]
		for any $\delta > 0$,
		\[
		\alpha > \frac{1}{2},\ p \geq \frac{4 \gamma}{\gamma - 1},
		\]
		\[
		\tvr(0, \cdot) = \vr(0, \cdot), \ \tvr \vU (0, \cdot) = \vm(0, \cdot).
		\]
		\item
		\[
		0 < \underline{\vr} \leq \tvr \leq \Ov{\vr} ,\ |\vU | \leq \Ov{\vu}
		\]	
		a.a. in $[0,T) \times \Td$.
		\item There exists $D \in L^1(0,T)$ such that
		\[
		\intTd{ \Big[ - \xi \cdot \vU (\tau, \cdot) (\xi \cdot \Grad ) \varphi + D(\tau) |\xi|^2 \varphi \Big] } \geq 0
		\ \mbox{for a.a.}\ \tau \in (0,T),
		\]
		for any $\xi \in R^d$ and any $\varphi\in C(\Td)$, $\varphi \geq 0$.	
	\end{itemize}
	Then
	\[
	\vr = \tvr,\ \vm = \tvr \vU \ \mbox{a.a. in}\ [0,T) \times \Td.
	\]		
	
\end{Proposition}

In particular, the (DW) solutions coincide with the smooth solutions introduced in Definition \ref{PE1}
emanating from the same initial data as long as the strong solution exists, meaning in the time interval $[0, T_{\rm max})$.

\begin{Remark} \label{rare}
The class of functions in which weak--strong uniqueness holds includes the planar rarefaction wave solutions, cf. \cite{FeiGhoJan}.
	\end{Remark}

\subsection{Euler system with uncertain data}

The main objective of the present paper is to develop a theoretical framework to deal with the Euler system \eqref{E1}--\eqref{E3} endowed with uncertain (random) data and its numerical approximations. In view of the generic
non-uniqueness of global in time solutions, it is convenient to work with solution \emph{sets} rather than individual solutions. Given a class of initial data $\mathcal{D}$, we consider the set $\mathcal{U}(\vr_0, \vm_0)$ of
all (DW) solutions emanating from $(\vr_0, \vm_0)$. The analysis is then performed on the set value mappings
\[
(\vr_0, \vm_0) \in \mathcal{D} \mapsto \mathcal{U}(\vr_0, \vm_0) \subset  \mathcal{T},
\]
where $\mathcal{T}$ is a suitable trajectory space.

Identifying suitable topologies on the data space, we reformulate some standard statistical tools, in
particular the Strong Law of Large Numbers (SLLN) in terms of solution sets. To this end, some recent
versions of SLLN for set--valued mappings ranging in Banach spaces will be used, see Ter{\' a}n \cite{Ter}
and the references cited therein.

Next, we consider a family of \emph{approximate} solutions $\vU_h (\vr_0, \vm_0)$, where $h \searrow 0$ is a discretization parameter. We show convergence of the approximate solutions towards the set
$\mathcal{U}(\vr_0, \vm_0)$ provided the approximation is stable consistent. Combining the abstract version
of SLLN with the bounds on the discretization error, we finally show convergence of the Monte Carlo method
for the isentropic Euler system in the class of (DW) solutions.

Convergence of the Monte Carlo approximations is {\it a priori} only weak due to possible and probably inevitable
oscillations developed in the approximation sequence. Fortunately, the
method of $\mathcal{K}-$convergence proposed in \cite{FeiLMMiz, FLSW} can be adapted to the random setting to
deduce strong convergence of the empirical averages of approximate solutions. Such a procedure can be seen as
``(Monte Carlo$)^2$'' method applied both at the level of exact statistical solutions and
their consistent approximations. Although the $\mathcal{K}-$convergence of consistent approximations and
the Monte Carlo method are based on the same idea of \emph{averaging},
the convergence of both is of different origin. The Monte Carlo method generates large families of \emph{i.i.d. samples} of random data,
for which the convergence of empirical averages follows from the Strong law of large numbers. The $\mathcal{K}-$convergence yields
convergence of numerical approximations on condition they are
asymptotically \emph{stationary} in the spirit of Birkhoff--Khinchin
ergodic theorem.

Finally, we address the problem of convergence towards the (unique) classical solution on its life span $[0, T_{\rm max})$. The abstract results are then applied to a finite volume approximation of the Euler system
proposed in \cite[Chapter 12]{FeLMMiSh} and the weak and strong convergence of the Monte Carlo finite volume method will be illustrated by numerical simulations.

To the best of our knowledge, this is the first rigorous analysis of convergence of the Monte Carlo
method for the multidimensional isentropic Euler system. We also refer to our recent works \cite{FLSY2022, LSY2024} where the convergence analysis of the Monte Carlo finite volume method for compressible Navier--Stokes(--Fourier) equations in the framework of global-in-time strong statistical solutions was presented.
The results for the isentropic Euler system lean on a combination
of several new ideas:

\begin{itemize}
	\item
	
	The concept of dissipative weak solutions as a universal closure of families of consistent approximations.
	
	\item Set--valued version of the Strong law of large numbers for general multivalued mapping with closed
	range.
	
	\item Application of Koml\' os version of convergence of empirical averages of integrable functions.
	
	\item Unconditionally convergent finite volume approximation scheme.
	
	\end{itemize}
	
	The paper is organized as follows. In Section \ref{D}, we introduce the topologies on both the data space and the trajectory space. Section \ref{S} reviews the properties
	of the multivalued solution set, in particular the closed graph property necessary for proving strong/weak measurability of the solution mapping, cf. Theorem \ref{TD1}. The problems with random data are introduced in Section~\ref{R}. Abstract (set--valued) version of SLLN is applied to the solution sets of the Euler system to obtain a general statement on convergence of the empirical averages of exact solutions emanating from i.i.d. data samples, see Theorem \ref{TR1}. In addition, its more standard (single--valued) version is proved on the maximal existence time of smooth solutions starting from smooth initial data, see Theorem \ref{R2a}. In Section \ref{C}, we consider consistent approximations of exact solutions. We show convergence of the Monte Carlo method for consistent approximations both in the weak form - Theorem \ref{TC1} --
	and the strong form -- Theorem  \ref{TC1s}.  We reformulate the previous results in the context of local smooth solutions in Section~\ref{cs}.
Finally, in Section~\ref{FV} we provide an example of a fully discrete unconditionally convergent finite volume scheme generating consistent approximations of the Euler system and illustrate theoretical results by numerical simulations for a well-known Kelvin-Helmholtz problem.

\section{Data and trajectory spaces}
\label{D}

In view of \eqref{E13}, a suitable data space can be chosen as
\begin{equation} \label{D1}
\mathcal{D} = \left\{ (\vr_0, \vm_0)\ \Big| \ (\vr_0, \vm_0)
\ \mbox{measurable on}\ \Td,\
\intTd{ E(\vr_0, \vm_0) } < \infty \right\}.
\end{equation}
Note that
\[
\mathcal{D} \subset L^\gamma(\Td) \times L^{\frac{2 \gamma}{\gamma + 1}}(\Td; R^d),
\]
and
\[
\vr_0(x) \geq 0 ,\ \mbox{and}\
\vr_0(x) = 0 \ \Rightarrow\ \vm_0 (x) = 0 \ \mbox{for a.a.}\ x \in \Td.
\]

Similarly, for a given $0 < T < \infty$, the trajectory space is
\begin{equation} \label{D2}
\mathcal{T} = \left\{ (\vr, \vm) \ \Big| \
(\vr, \vm) \in C_{\rm weak}([0,T]; L^\gamma (\Td) \times L^{\frac{2 \gamma}{\gamma + 1}}(\Td; R^d)) \right\}.
	\end{equation}

The following result was shown in \cite[Section 3.1]{BreFeiHof19}.

\begin{Proposition}[{\bf Weak sequential stability}] \label{PS1}
	
	Suppose $(\vr_{0,n}, \vm_{0,n} )_{n=1}^\infty \subset \mathcal{D}$ is a sequence of initial data
	satisfying
	\begin{align}
	\vr_{0,n} &\to \vr_0 \ \mbox{in}\ W^{-\ell,2}(\Td)),\
	\vm_{0,n} \to \vm_0 \ \mbox{in}\ W^{-\ell,2}(\Td; R^d),\br
	\intTd{ E(\vr_{0,n}, \vm_{0,n}) } &\to \intTd{ E(\vr_0, \vm_0 ) }.
	\label{D4}
	\end{align}
	Let $(\vr_n, \vm_n)_{n=1}^\infty$,
	\[
	(\vr_n, \vm_n) \in \mathcal{U}(\vr_{0,n}, \vm_{0,n})
	\]
	be a sequence of the corresponding (DW) solutions of the Euler system.
	
	Then $(\vr_0, \vm_0) \in \mathcal{D}$, and
	there is a subsequence $(\vr_{n_k}, \vm_{n_k})_{k=1}^\infty$,
	\[
	\vr_{n_k} \to \vr \ \mbox{in} \ C_{\rm weak}([0,T]; L^\gamma(\Td)),\
	\vm_{n_k} \to \vm \ \mbox{in}\ C_{\rm weak}([0,T]; L^{\frac{2 \gamma}{\gamma + 1}}(\Td; R^d)) \ \mbox{as}\ k \to \infty,
	\]
	where $(\vr, \vm) \in \mathcal{U}(\vr_0, \vm_0)$.

\end{Proposition}

\subsection{Metrics on the data space}

We define a metrics on the data space,
\begin{equation} \label{D5}
d_{\mathcal{D}}	\Big[ (\vr^1_0, \vm^1_0); (\vr^2_0, \vm^2_0) \Big] =
\| \vr^1_0 - \vr^2_0 \|_{L^\gamma(\Td)} + \left\| \mathds{1}_{\vr^1_0 > 0} \frac{\vm^1_0}{\sqrt{\vr^1_0}} -
\mathds{1}_{\vr^2_0 > 0} \frac{\vm^2_0}{\sqrt{\vr^2_0}}  \right\|_{L^2(\Td; R^d)}.
	\end{equation}
It is easy to check that convergence in the metrics $d_{\mathcal{D}}$ is equivalent to the convergence stated
in \eqref{D4}, and that $[ \mathcal{D}, d_{\mathcal{D}} ]$ is a (metric) Polish space.

\subsection{Topology of the trajectory space}
\label{ssc}

A vast majority of available results on set--valued SLLN require
the trajectory space to be endowed with a topology of a separable Banach space. We distinguish the \emph{weak} setting, where
$\mathcal{T}$ is viewed as a subspace of the Banach space
\[
\mathcal{T}_{\rm weak} \hookrightarrow C([0,T]; W^{-\ell,2}(\Td; R^{d+1})), \ \ell > d,
\]
and the \emph{strong} setting, where
\[
\mathcal{T}_{\rm strong} \hookrightarrow L^q((0,T) \times \Td; R^{d+1}),\ 1 < q \leq \frac{2 \gamma}{\gamma + 1}
= \min \left\{ \gamma; \frac{2 \gamma}{\gamma + 1} \right\}.
\]

\section{Properties of solution sets}
\label{S}
Given initial data $(\vr_0, \vm_0) \in \mathcal{D}$, we review the basic properties of the associated \emph{solutions set}:
\begin{align}
	\mathcal{U}(\vr_0, \vm_0)= \Big\{ (\vr, \vm) \in \mathcal{T}\ \Big| \
	&(\vr, \vm) \ \mbox{a (DW) solution of the Euler system in}\  [0,T] \times \Td ,\br &(\vr, \vm)(0, \cdot) = (\vr_0, \vm_0) \Big\}. \label{S1}
\end{align}	

\subsection{Topological properties of the solution set}

For each $(\vr_0, \vm_0) \in \mathcal{D}$, there holds:
\begin{itemize}
	\item The set $\mathcal{U}(\vr_0, \vm_0)$ is non--empty, see Proposition \ref{PE4}.
	\item The set $\mathcal{U}(\vr_0, \vm_0)$ is convex.
	\item The set $\mathcal{U}(\vr_0, \vm_0)$ is a closed bounded subset of $L^\infty(0,T; L^\gamma(\Td) \times L^{\frac{2 \gamma}{\gamma + 1}}(\Td; R^d))$.
	\item The set $\mathcal{U}(\vr_0, \vm_0)$ is a compact subset of $C([0,T]; W^{-\ell,2}(\Td) \times W^{-\ell,2}(\Td; R^d))$, $\ell > d$.
	\item The set $\mathcal{U}(\vr_0, \vm_0)$ is a closed bounded subset of $L^q((0,T) \times \Td; R^{d+1})$,
	$1 \leq q \leq \frac{2\gamma}{\gamma + 1}$.

\end{itemize}

\subsection{Measurability of the solution mapping in the
weak topology}
\label{MDM}

By virtue of weak sequential stability established in Proposition \ref{PS1}, the solution mapping
\[
(\vr_0, \vm_0) \in \mathcal{D} \mapsto \mathcal{U}(\vr_0, \vm_0) \subset C([0,T]; W^{-\ell,2}(\Td; R^{d+1})),\ \ell > d,
\]
enjoys the closed graph property. Applying \cite[Lemma 12.1.8]{StrVar} we conclude that
the solution mapping is (strongly) Borel measurable with respect
to the topologies of $[\mathcal{D}; d_{\mathcal{D}}]$ and of the Banach space
$C([0,T]; W^{-\ell,2}(\Td; R^{d+1}))$. Equivalently, we can say that the set--valued mapping
\[
(\vr_0, \vm_0) \in [\mathcal{D}; d_{\mathcal{D}}] \mapsto \mathcal{U}(\vr_0, \vm_0)
\in {\rm comp}[ C([0,T]; W^{-\ell,2}(\Td; R^{d+1})]
\]
ranging in the set
${\rm comp}[ C([0,T]; W^{-\ell,2}(\Td; R^{d+1})]$ of all compact subsets of the Banach space
$C([0,T]; W^{-\ell,2}(\Td; R^{d+1}))$ is Borel measurable with respect to the Hausdorff topology on the space ${\rm comp}[ C([0,T]; W^{-\ell,2}(\Td; R^{d+1})]$.

\subsection{Measurability in the strong $L^q$-topology}
\label{MDMS}

Next, we address the problem of measurability of the solution mapping if the target trajectory space is endowed with the (strong) topology
of the separable reflexive Banach space $L^q((0,T) \times \Td; R^{d+1})$, $1 < q \leq \frac{2 \gamma}{\gamma + 1}$.

Suppose the data $(\vr_0, \vm_0)$ is a random variable defined on a complete probability space $(\Omega, \mathfrak{M}, \prst)$ ranging in $\mathcal{D}$. In accordance with Section \ref{MDM},
the mapping $\mathcal{U} \in {\rm comp}[ C([0,T]; W^{-\ell,2}(\Td; R^{d+1}))]$ is a set--valued random variable, meaning a random closed set.
Let $(\bfphi_i)_{i = 1}$ be a Castaing representation of $\mathcal{U}$, meaning a countable family of selections (random variables)
\begin{align}
	\bfphi_i: \Omega &\to  C([0,T]; W^{-\ell,2}(\Td; R^{d+1})),\ \bfphi_i \in \mathcal{U},\br
	{\bf cl}_{C([0,T]; W^{-\ell,2}(\Td; R^{d+1})}[(\bfphi_i)_{i=1}^\infty] &= \mathcal{U}.
	\label{S2}
\end{align}
As the sets $\mathcal{U}$ are bounded in $L^q((0,T) \times \Td; R^{d + 1})$ $\prst$-a.s., we get
\[
{\bf cl}_{C([0,T]; W^{-\ell,2}(\Td; R^{d+1})}[(\bfphi_i)_{i=1}^\infty] =
{\bf cl}_{L^q((0,T) \times \Td; R^{d + 1})_{\rm weak}}[(\bfphi_i)_{i=1}^\infty] = \mathcal{U}
\]
Moreover, since the sets $\mathcal{U}$ are convex, we may introduce a new countable family of random variables,
\[
\bfomega_{\lambda} = \sum_{i=1}^N \lambda_i \bfphi_i, \
\lambda = (\lambda_1,\dots, \lambda_N) \in Q^N, \ \lambda_i \geq 0,\ \sum_{i=1}^N \lambda_i = 1,\ N = 1,2,\dots
\]
where $\bfomega_\lambda \in \mathcal{U}$. By virtue of Banach-Saks theorem,
\[
{\bf cl}_{L^q((0,T) \times \Td; R^{d + 1})}[(\bfomega_\lambda)_\lambda] = {\bf cl}_{L^q((0,T) \times \Td; R^{d + 1})_{\rm weak}}[(\bfphi_i)_{i=1}^\infty] = \mathcal{U}.
\]
In other words, the family $(\bfomega_\lambda)_\lambda$ is a Castaing representation of $\mathcal{U}$ in the (strong) topology of $L^q((0,T) \times \Td; R^{d + 1})$. We conclude the set mapping
\[
\omega \in \Omega \mapsto (\vr_0, \vm_0) \in \mathcal{D} \mapsto \mathcal{U}(\vr_0, \vm_0)
\]
is (strongly) measurable with respect to the (strong) topology of $L^q((0,T) \times \Td; R^{d + 1})$, $1 < q \leq \frac{2 \gamma}{\gamma + 1}$.

\section{Euler system with random initial data}
\label{R}

We consider the Euler system with the initial data $(\vr_0, \vm_0)$ being random on a complete probability basis $(\Omega, \mathfrak{M}, \prst)$. We focus on the set--valued solution mapping
\[
\omega \in \Omega \mapsto (\vr_0, \vm_0) (\omega) \in \mathcal{D} \mapsto \mathcal{U}(\vr_0, \vm_0)(\omega)	\in 2^{\mathcal{T}},
\]
where $2^{\mathcal{T}}$ denotes the family of all subsets of the trajectory space $\mathcal{T}$.

In accordance with Section \ref{ssc}, we consider two topologies on the trajectory space $\mathcal{T}$:
\begin{itemize}
	\item The space $\mathcal{T}_{\rm weak}$ with the topology of the Banach space
	\[
	C([0,T]; W^{-\ell,2}(\Td; R^d)).
	\]
	The subspace of all compact sets
	\[
	{\rm comp}[\mathcal{T}_{\rm weak}] \subset 2^\mathcal{T}
	\]
	is converted to a Polish space with the Hausdorff metrics
	\[
	D_{\mathcal{H}}[ \mathcal{A}; \mathcal{B} ] = \max \left[ \max_{y \in \mathcal{A}} {\rm dist}_{\mathcal{T}_{\rm weak}}[y, \mathcal{B} ] ; \max_{y \in \mathcal{B}} {\rm dist}_{\mathcal{T}_{\rm weak}}[y, \mathcal{A} ]   \right]
	\]
	for $\mathcal{A}, \mathcal{B} \in {\rm comp}[\mathcal{T}_{\rm weak}]$.
	\item The space $\mathcal{T}_{\rm strong}$ with the topology of the Banach space
	\[
	L^q((0,T) \times \Td; R^{d+1}),\ 1 < q \leq \frac{2 \gamma}{\gamma + 1}.
	\]
	The space of all closed subsets
	\[
	{\rm closed}[\mathcal{T}_{\rm strong}] \subset 2^\mathcal{T}
	\]
	is endowed with the Wijsman topology, with the convergence
	\[
	\mathcal{A}_n \to \mathcal{A} \ \Leftrightarrow \ {\rm dist}_{\mathcal{T}_{\rm strong}} [ y, \mathcal{A}_n]
	\to {\rm dist}_{\mathcal{T}_{\rm strong}} [ y, \mathcal{A}] \ \mbox{for any}\ y \in \mathcal{T}.
	\]
	\end{itemize}

By Hess' measurability theorem \cite{Hess}, (strong) measurability of a closed set--valued mapping
is equivalent to its Borel measurability with respect to the Wijsman topology. Consequently,
summing up the material of Sections \ref{MDM}, \ref{MDMS} we obtain the following result.

	\begin{Theorem}[{\bf Measurability of the solution sets}] \label{TD1} Suppose the initial data
		$(\vr_0, \vm_0)$ is a random variable ranging in the Polish space $[ \mathcal{D}; d_{\mathcal{D}}]$.
		
		Then the solution mapping
		\[
		\omega \in \Omega \mapsto (\vr_0, \vm_0) \in \mathcal{D} \mapsto \mathcal{U} (\vr_0, \vm_0) \in 2^{\mathcal{T}}
		\]
		is measurable as the set--valued mapping ranging in ${\rm comp}[\mathcal{T}_{\rm weak}]$ endowed with
		the Hausdorff topology or in ${\rm closed}[\mathcal{T}_{\rm strong}]$ endowed with the Wijsman topology.

		\end{Theorem}

\subsection{Strong law of large numbers}

With Theorem \ref{TD1} at hand, we may apply the available results concerning SLLN for set--valued mappings.
Consider a sequence of pairwise independent equally distributed (i.i.d.) random data
\begin{equation} \label{R1}
	(\vr_{0,n}, \vm_{0,n}) \in \mathcal{D} \ \mbox{with the associated solution sets}\ \mathcal{U}_n = \mathcal{U}(\vr_{0,n}, \vm_{0,n}) \in 2^{\mathcal{T}}.
\end{equation}
For a random set $\mathcal{U}$, its Aumann expectation is defined as
\[
\expe{ \mathcal{U} } = {\rm cl} \left\{ \expe{ \xi } \ \Big|\ \xi \ \mbox{an integrable selection of} \ \mathcal{U} \right\}.
\]
Here and hereafter, the symbol $\mathbb{E}$ denotes expectation.

\begin{Theorem}\label{TR1} {\bf (Strong law of large numbers for random Euler system).}
	
	\noindent
	Let $(\vr_{0,n}, \vm_{0,n})_{n=1}^\infty$ be a sequence of pairwise i.i.d. copies of  random data $(\vr_0,
	\vm_0) \in \mathcal{D}$, with the associated sequence of
	sets of (DW) solutions $\left( \mathcal{U}(\vr_{0,n}, \vm_{0,n} ) \right)_{n=1}^\infty$. Suppose
	\[
	\expe{ \intTd{ E(\vr_0, \vm_0) } } < \infty.
	\]
	
	Then
	\begin{itemize}
		\item
.		\begin{equation} \label{R2}
			\frac{1}{N} \sum_{n=1}^N \mathcal{U}(\vr_{0,n}, \vm_{0,n}) \to \expe{ \mathcal{U}(\vr_{0}, \vm_{0}) } \
			\mbox{as}\ N \to \infty\  \mbox{a.s.}
		\end{equation}	
		in the Hausdorff topology on compact subsets of $C([0,T]; W^{-\ell,2}(\Td; R^{d+1}))$;
		\item
		\begin{equation} \label{R3}
			\frac{1}{N} \sum_{n=1}^N \mathcal{U}(\vr_{0,n}, \vm_{0,n}) \to \expe{ \mathcal{U}(\vr_{0}, \vm_{0}) } \
			\mbox{as}\ N \to \infty \  \mbox{a.s.}
		\end{equation}	
		in the Wijsman topology on closed convex subsets of the Lebesgue space $L^q((0,T) \times \Td; R^{d+1})$, $1 < q \leq \frac{2 \gamma}{\gamma + 1} $. 	
	\end{itemize}	
\end{Theorem}

The statement \eqref{R2} can be found in Artstein and Hansen \cite{ArtHan}, its generalization to non--compact
closed sets stated in \eqref{R3} was proved by Ter\' an \cite{Ter}. As a matter of fact, Ter\' an's version
of SLLN holds for general closed not necessarily convex set--valued mappings with respect a stronger \emph{gap} topology. For convex set--valued mappings, however, the gap and Wijsman topologies (with respect to arbitrary equivalent
norm) coincide.

\subsection{Convergence for smooth data}

In addition to \eqref{E14}, suppose the initial data $(\vr_0, \vm_0)$ are smooth as in Proposition \ref{PE1},
\begin{equation} \label{E5a}
	\vr_0 \in W^{k,2}(\Td),\ \inf_{x \in \Td} \vr_0 (x) > 0,\ \vm_0 \in W^{k,2}(\Td; R^d),\ k > 1 + \frac{d}{2}.
\end{equation}	
As stated in Proposition \ref{PE1}, the Euler system admits a classical solution defined on a maximal time interval
$[0, T_{\rm max})$, $T_{\rm max} \leq \infty$.
We introduce the data space for strong data,
\begin{equation} \label{E11}
	\mathcal{D}_S = \left\{ (\vr_0, \vm_0) \ \Big| \ \vr_0 \in W^{k,2}(\Td),\ \inf_{x \in \Td} \vr_0(x) > 0,\
	\vm_0 \in W^{k,2}(\Td; R^d) \right\} \subset \mathcal{D},
	\end{equation}
with a metrics
\begin{equation} \label{E12}
d_{\mathcal{D}_S} \left[ (\vr^1_0, \vm^1_0); (\vr^2_0, \vm^2_0) \right] = \| \vr^1_0 - \vr^2_0 \|_{W^{k,2}(\Td)}  +
\| \vm^1_0 - \vm^2_0 \|_{W^{k,2}(\Td;R^d)} + \left\| \frac{1}{\vr^1_0} - \frac{1}{\vr^2_0} \right\|_{C(\Td)}.
	\end{equation}
Following \cite[Section 2.1]{FeiLuk2022}
we can show lower semi--continuity of the
blow up time $T_{\rm max}(\vr, \vm)$ of a strong solution $(\vr, \vm)$ exactly as in \cite[Theorem 2.1]{FeiLuk2022}.

\begin{Proposition}[{\bf Lower semi--continuity of $T_{\rm max}$}] \label{CE2}
	The mapping
	\[
	(\vr_0, \vm_0) \in \mathcal{D}_S \mapsto T_{\rm max}(\vr, \vm) \in (0, \infty)
	\]
	is lower semi--continuous.
	
	\end{Proposition}

\subsubsection{SLLN for smooth data}

Suppose the data $(\vr_0, \vm_0)$ are random and belong to the class \eqref{E5a} a.s. As a consequence of
lower semi--continuity of $T_{\rm max}$ stated in Proposition \ref{CE2}, the cut--off solutions
\[
\mathds{1}_{t \leq T_{\rm max}} (\vr, \vm)(t,x),\ t \in [0, \infty),\ x \in \Td
\]
are random variables. Applying the standard version of SLLN on separable Banach spaces (see
Etemadi \cite{Etem}  ) we obtain the following result:

\begin{Theorem}\label{TR1a} {\bf (Strong law of large numbers for random Euler system, strong solutions.)}
	
	\noindent
	Let $(\vr_{0,n}, \vm_{0,n})_{n=1}^\infty$ be a sequence of pairwise i.i.d. representations of smooth random data $(\vr_0,
	\vm_0) \in \mathcal{D}_S$, with the associated sequence of
	smooth solutions $(\vr_n, \vm_n )_{n=1}^\infty$. Suppose
	\[
	\expe{ \intTd{ E(\vr_0, \vm_0) } } < \infty.
	\]
	
	Then
		\begin{equation} \label{R2a}
			\frac{1}{N} \sum_{n=1}^N \mathds{1}_{t \leq T_{\rm max}}(\vr_{n}, \vm_{n}) \to \expe{\mathds{1}_{t \leq T_{\rm max}} (\vr, \vm) } \ \ \mbox{as}\ N \to \infty \  \mbox{a.s.}
		\end{equation}	
		in $C([0,T]; W^{-\ell,2}(\Td; R^{d+1}))$ and  $L^q((0,T) \times \Td; R^{d+1})$, $1 < q \leq \frac{2 \gamma}{\gamma + 1} $, where $(\vr, \vm)$ is the unique smooth solution emanating from the data
		$(\vr_0, \vm_0)$. 	
	
\end{Theorem}

\begin{Remark} \label{Rcs1}
	The maximal existence time $T_{\rm max}$ is bounded below by the norm of the initial data,
	\[
	T_{\rm max} \geq C \left( \left\| \vr_0^{\frac{\gamma - 1}{2}} \right\|_{W^{k,2}(\Td)}  +
	\| \vu_0 \|_{W^{k,2}(\Td; R^d)} \right)^{-2},\ \vu_0 = \frac{\vm_0}{\vr_0},
	\]
	see e.g. \cite[Theorem 5.1]{GuoWu}. In particular, $T_{\rm max}$ is bounded below away from zero
	by a positive \emph{deterministic} constant if the norm of the data
	\[
	\left\| \vr_0^{\frac{\gamma - 1}{2}} \right\|_{W^{k,2}(\Td)}  +
	\| \vu_0 \|_{W^{k,2}(\Td; R^d)}
	\] 	
is bounded above by a deterministic constant.	
\end{Remark}

\section{Consistent approximations}
\label{C}

In the preceding part, we have derived several forms of SLLN for
\emph{exact} (DW) solutions of the isentropic Euler system. In this section, we focus on consistent approximations.

\begin{Definition}[\bf Consistent approximation]\label{D_CA}
Suppose the initial data $(\vr_0, \vm_0) \in \mathcal{D}$
satisfy
\begin{equation} \label{C1a}
\frac{1}{L} \leq \vr_0(x) \leq L,\ |\vm_0(x) | \leq L \ \mbox{for a.a.}\ x \in \Td
\end{equation}
for some constant $L > 0$. 	
The consistent approximation of the Euler system with the initial data $(\vr_0, \vm_0)$ is a family $(\vr_h, \vm_h)_{h \searrow 0}$ satisfying
\begin{equation} \label{C1}
(\vr_h, \vm_h) \in L^\infty((0,T) \times \Td; R^{d+1}),\ \inf_{(0,T) \times \Td} \vr_h > 0;
\end{equation}	
\begin{equation} \label{C2}
	\int_0^T \intTd{ \Big[ \vr_h \partial_t \varphi + \vm_h \cdot \Grad \varphi \Big] } \dt =  - \intTd{ \vr_0 \varphi(0, \cdot)} + e^1(h,L, \| \varphi \|_{C_2} )
\end{equation}
for any $\varphi \in C^2_c([0,T) \times \Td)$;
	\begin{align}
	\int_0^T &\intTd{ \left[ \vm_h \cdot \partial_t \bfphi + \frac{\vm_h \otimes \vm_h}{\vr_h} : \Grad \bfphi +
		a \vr_h^\gamma \Div \bfphi \right] } \dt \br &= - \intTd{ \vm_0 \cdot \bfphi(0, \cdot) } + e^2(h,L, \| \bfphi \|_{C^2})
	\label{C3}
\end{align}
for any $\bfphi \in C^2_c([0,T) \times \Td; R^d)$;
\begin{equation} \label{C4}
	\intTd{ E(\vr_h, \vm_h) (\tau, \cdot) } \leq \intTd{ E(\vr_0, \vm_0) } + e^3(h,L)
	\end{equation}
for a.a. $\tau \in (0,T)$.
The consistency error terms satisfy
\begin{equation} \label{C5}
e^1(h, L, \| \varphi \|_{C^2}) \to 0,\
e^2(h, L, \| \bfphi \|_{C^2}) \to 0,\ e^3(h, L) \to 0 \ \mbox{as}\ h \to 0
	\end{equation}
for any fixed $L$ and $\| \varphi \|_{C^2}$, $\| \bfphi \|_{C^2}$.
\end{Definition}

 As presented in \cite{FeLMMiSh} consistent approximations are generated by suitable structure preserving \emph{numerical schemes}.
 Here and hereafter, we suppose the consistent approximation is determined by a mapping $\vU_h$,
\[
(\vr_h, \vm_h) = \vU_h (\vr_0, \vm_0),\ \mbox{where}\  \vU_h : \mathcal{D} \to L^1((0,T) \times \Td; R^{d+1})
\ \mbox{is Borel measurable.}
\]

\subsection{Convergence of consistent approximations}

The following result was proved in \cite[Chapter 7, Theorem 7.9]{FeLMMiSh}.

\begin{Proposition}[{\bf Convergence of consistent approximations}] \label{PC1}
Let the initial data $(\vr_0, \vm_0) \in \mathcal{D}$ satisfying \eqref{C1a} be given. Let
$(\vU_h)_{h \searrow 0}$ be a consistent approximation.

Then any sequence $h_n \to 0$ contains a subsequence $(h_{n_k})_{k=1}^\infty$ such that
\begin{align} \
\vU_{h_{n_k}} (\vr_0, \vm_0) &\to (\vr, \vm) \ \mbox{weakly-(*) in}\
L^\infty(0,T; L^\gamma(\Td) \times L^{\frac{2 \gamma}{\gamma + 1}}(\Td; R^d)), \br \mbox{where}\
(\vr, \vm) &\in \mathcal{U}(\vr_0, \vm_0).
\label{C6}	
	\end{align}
	\end{Proposition}

As the weak-(*) topology is metrizable on bounded sets of $L^\infty(0,T; L^\gamma(\Td) \times L^{\frac{2 \gamma}{\gamma + 1}}(\Td; R^d))$, we deduce from Proposition \ref{PC1} that
\begin{equation} \label{C6a}
{\rm dist}_{L^\infty(0,T; L^\gamma(\Td) \times L^{\frac{2 \gamma}{\gamma + 1}}(\Td; R^d))- ({\rm weak-(*))}}	
\left[ \vU_h (\vr_0, \vm_0), \mathcal{U}(\vr_0, \vm_0) \right] \to 0
\ \mbox{as}\ h \to 0
	\end{equation}
whenever $(\vr_0, \vm_0) \in \mathcal{D}$ satisfy \eqref{C1a} and $\vU_h$ is a consistent approximation.

\subsection{Monte Carlo method, weak form}

Suppose that the data $(\vr_0, \vm_0) \in \mathcal{D}$ are random. For the sake of simplicity, we assume
\begin{equation} \label{C1b}
	\frac{1}{L} \leq \vr_0(x) \leq L,\ |\vm_0(x) | \leq L \ \mbox{for a.a.}\ x \in \Td
\end{equation}
$\prst-$a.s., where $L$ is a \emph{deterministic} constant. In particular
\begin{equation} \label{C1c}
	\expe{ \intTd{ E(\vr_0, \vm_0 ) }} < \infty.
\end{equation}	

Consider $(\vr_{0,n}, \vm_{0,n})_{n=1}^\infty$ -  a sequence of pairwise independent copies of $(\vr_0, \vm_0)$ -
along with a consistent approximation $\vU_h$. As $\vU_h$ is Borel measurable, the approximations
$(\vU_h (\vr_{0,n}, \vm_{0,n}))_{n=1}^\infty$ are pairwise independent equally distributed,
\[
\vU_h (\vr_{0,n}, \vm_{0,n}) \sim \vU_h (\vr_0, \vm_0),
\]
where the symbol $\sim$ stands for equivalence in law.

Our goal is to estimate the distance between the Monte Carlo estimator using consistent approximations
\[
\frac{1}{N} \sum_{n=1}^N \vU_h (\vr_{0,n}, \vm_{0,n})
\]
and the expected value of the solution set
\[
\expe{ \mathcal{U}(\vr_0, \vm_0)}.
\]
Note carefully that the approximate solutions
\[
\frac{1}{N} \sum_{n=1}^N \vU_h (\vr_{0,n}, \vm_{0,n})
\]
are single--valued random variables ranging in the trajectory space, while the limit $\expe{ \mathcal{U}(\vr_0, \vm_0)}$ is a (convex)
subset of the same space.

Write, first formally,
\begin{align}
{\rm dist} \left[ \frac{1}{N} \sum_{n=1}^N \vU_h (\vr_{0,n}, \vm_{0,n}) ;
\expe{ \mathcal{U}(\vr_0, \vm_0)} \right] &\leq
{\rm dist} \left[ \frac{1}{N} \sum_{n=1}^N \vU_h (\vr_{0,n}, \vm_{0,n}) ;
\frac{1}{N} \sum_{n=1}^N \mathcal{U} (\vr_{0,n}, \vm_{0,n}) \right] \br
&+ {\rm dist} \left[\frac{1}{N} \sum_{n=1}^N \mathcal{U} (\vr_{0,n}, \vm_{0,n});\expe{ \mathcal{U}(\vr_0, \vm_0)} \right],
\label{CC1}
\end{align}
where the former term on the right--hand side represents the
\emph{discretization} error, while the latter is the \emph{statistical} error.

The discretization error is controlled by \eqref{C6a}. Moreover,
as the data are deterministically bounded by \eqref{C1a} and
obey the uniform stability bounds stated in \eqref{C4}, we may infer
\[
\expe{{\rm dist}_{L^\infty(0,T; L^\gamma(\Td) \times L^{\frac{2 \gamma}{\gamma + 1}}(\Td; R^d))- ({\rm weak-(*))}}	
	\left[ \vU_h (\vr_0, \vm_0), \mathcal{U}(\vr_0, \vm_0) \right] }
	\to 0 \ \mbox{as}\ h \to 0.
\]
Consequently, as all data are equally distributed we conclude
\begin{align}
&\expe{{\rm dist} \left[ \frac{1}{N} \sum_{n=1}^N \vU_h (\vr_{0,n}, \vm_{0,n}) ;
\frac{1}{N} \sum_{n=1}^N \mathcal{U} (\vr_{0,n}, \vm_{0,n}) \right]} \br &\leq \frac{1}{N} \sum_{n=1}^N \expe{ {\rm dist} \Big[ \vU_h (\vr_{0,n}, \vm_{0,n}); \mathcal{U} (\vr_{0,n}, \vm_{0,n})  \Big] } \br
&= \expe{{\rm dist}_{L^\infty(0,T; L^\gamma(\Td) \times L^{\frac{2 \gamma}{\gamma + 1}}(\Td; R^d))- ({\rm weak-(*))}}	
	\left[ \vU_h (\vr_0, \vm_0), \mathcal{U}(\vr_0, \vm_0) \right] }
\to 0 \ \mbox{as}\ h \to 0.
\label{CC2}
\end{align}

The convergence of
the statistical error on the right-hand side of \eqref{CC1} follows from Theorem \ref{TR1}. Adjusting
the metrics on the trajectory space, we get the following result.

\begin{Theorem}\label{TC1} {\bf (Convergence of Monte Carlo method, weak form).}
	
\noindent	
Suppose $(\vr_{0,n}, \vm_{0,n})_{n=1}^\infty$ is a family of pairwise i.i.d. copies
of initial data $(\vr_0, \vm_0) \in \mathcal{D}$ satisfying \eqref{C1a} with a deterministic constant $L$.
Let $(\vU_h)_{h \searrow 0}$ be a consistent approximation  in the sense of Definition~\ref{D_CA}.

Then
\begin{equation} \label{C7}
\expe{ {\rm dist}_{W^{-\ell,2}((0,T) \times \Td; R^{d+1}) } \left[	\frac{1}{N} \sum_{n=1}^N \vU_h (\vr_{0,n}, \vm_{0,n})  ; \expe{ \mathcal{U}(\vr_{0}, \vm_{0})   }\right] }
\to 0 \ \mbox{as}\ h \to 0,\ N \to \infty
	\end{equation} 	
for any $\ell > d + 1$.	
	\end{Theorem}

\subsection{Monte Carlo method, strong form}

The convergence result stated in Theorem \ref{TC1} holds in the weak topology of the trajectory space, whereas the limit is a set
rather than a single function.
In principal, this is optimal in view of the weak convergence of the consistent approximations and ``generic'' non--uniqueness of the weak solutions to the Euler system.

As shown in \cite[Chapter 7, Theorem 7.9]{FeLMMiSh}, the weak convergence can be converted to the strong one up to
a suitable subsequence
by replacing the consistent approximations by their
empirical averages in the spirit of the classical Banach--Saks theorem.

\begin{Proposition}[{\bf ($\mathcal{K}$)-convergence of consistent approximations}] \label{PC2}
	
Let the initial data \linebreak $(\vr_0, \vm_0) \in \mathcal{D}$ satisfying \eqref{C1a} be given. Let
$(\vU_h)_{h \searrow 0}$ be a consistent approximation.

Then any sequence $h_n \to 0$ contains a subsequence $(h_{n_k})_{k=1}^\infty$ such that 	
	\begin{equation} \label{C8}
\frac{1}{K} \sum_{k=1}^K \vU_{h_{n_k}}(\vr_0, \vm_0) \to (\vr, \vm) \in \mathcal{U}(\vr_0, \vm_0) \ \mbox{in}\
L^q((0,T) \times \Td; R^{d+1}) \ \mbox{as}\ K \to \infty
\end{equation}
for any $1 < q \leq \frac{2 \gamma}{\gamma + 1}$.	
Moreover, the convergence \eqref{C8} holds for any subsequence of $(h_{n_k})_{k=1}^\infty$.
\end{Proposition}

Pursuing the idea of Balder \cite{Bald, BALDER},
we may extend validity of \eqref{C8} to
	\begin{equation} \label{C8a}
	\frac{1}{K} \sum_{k=1}^K B( \vU_{h_{n_k}}(\vr_0, \vm_0) ) \to \Ov{B (\vr, \vm)},  \ \mbox{in}\
	L^1((0,T) \times \Td) \ \mbox{as}\ K \to \infty
\end{equation}
for any (globally) Lipschitz $B: R^{d+1} \to R$.
Relation \eqref{C8a} identifies an analogue of the standard Young measure
called $(S)$-limit in \cite{Fei2021SC}, see also \cite{FLSS, FeiLMMiz,  FLSW} for applications in numerical methods.

Suppose now that $(\vr_0, \vm_0)$ are random.
The crucial observation is that the $(S)$-limit identified in
\eqref{C8a} can be obtained for the same subsequence
\emph{independently} of the random event $\omega \in \Omega$.
This can be shown by means of
the celebrated refinement of Banach--Saks theorem due to Koml\' os \cite{Kom}.
Similarly to the preceding section, suppose the data $(\vr_0, \vm_0)$ satisfy \eqref{C1a} a.s. with a deterministic
constant $L$. Let $h_n \to 0$ be an arbitrary sequence of vanishing
discretization parameters. It follows from the approximate energy inequality
\eqref{C4} that
\[
\sup_n \expe{ \int_0^T \intTd{ E(\vr_{h_n}, \vm_{h,n}) }  \dt } =
\sup_n \int_{\Omega} \int_0^T \intTd{ E(\vr_{h_n}, \vm_{h,n}) } \dt \ \D \omega < \infty,
\]
where $(\vr_{h_n}, \vm_{h_n}) = \vU_{h_n}(\vr_0, \vm_0)$.
In particular,
\[
\sup_n \left\| B \left( \vU_{h_n}(\vr_0, \vm_0) \right) \right\|_{L^1(\Omega \times (0,T) \times \Td; R^{d+1})} < \infty
\]
for any globally Lipschitz $B: R^{d + 1} \to R$.

Consequently, by virtue of Koml\' os theorem, there exists a subsequence $(h_{n_k})_{k=1}^\infty$ such that
\[
\frac{1}{K} \sum_{k=1}^K B \left( \vU_{h_{n_k}} (\vr_0, \vm_0) \right) \to \Ov{B(\vr, \vm)} \ \mbox{a.a. in}\  \Omega \times (0,T) \times \Td, \
\mbox{and in}\ L^1(\Omega \times (0,T) \times \Td)
\]
for any Lipschitz $B$.
In particular, there is a measurable subset of random elements
$\Omega'$ of full measure such that
\begin{equation} \label{random}
\frac{1}{K} \sum_{k=1}^K B(\vU_{h_{n_k}} (\vr_0, \vm_0)(\omega)) \to \Ov{B(\vr, \vm)} \ \mbox{a.a. in}\
(0,T) \times \Td \ \mbox{for any}\ \omega \in \Omega',
\end{equation}
meaning the convergence in \eqref{random} takes place $\prst-$a.s.

Finally, for linear $B$, we get, in particular,
\begin{equation} \label{ran2}
\frac{1}{K} \sum_{k=1}^K \vU_{h_{n_k}} (\vr_0, \vm_0) \to (\vr, \vm) \ \mbox{in}\
L^q((0,T) \times \Td; R^{d+1}) \  \mbox{a.s.}
\end{equation}
for any $1 < q \leq \frac{2 \gamma}{\gamma + 1}$, where $(\vr, \vm) \in \mathcal{U}(\vr_0, \vm_0)$. The limit, being a pointwise limit of measurable functions, is measurable.

Summarizing the above discussion we deduce the following conclusion.

\begin{Proposition}\label{PC3} {\bf (($\mathcal{K}$)-convergence of consistent approximations, random setting).}
	
	\noindent
	Let $(\vr_0, \vm_0) \in \mathcal{D}$ be random data satisfying \eqref{C1a} with a deterministic
	constant $L$. Let $(\vU_h)_{h \searrow 0}$ be a consistent approximation of the Euler system.
	
	Then any sequence $h_n \to 0$ contains a subsequence $(h_{n_k})_{k=1}^\infty$ such that
	\begin{equation} \label{C12}
		\frac{1}{K} \sum_{k=1}^K \vU_{h_{n_k}}(\vr_0, \vm_0) \to (\vr, \vm) \in \mathcal{U}(\vr_0, \vm_0) \ \mbox{in}\
		L^q((0,T) \times \Td; R^{d+1}) \ \mbox{as}\ K \to \infty
	\end{equation}
	for any $1 \leq q < \frac{2 \gamma}{\gamma + 1}$, a.s. Moreover
	\begin{equation} \label{C12a}
	\frac{1}{K} \sum_{k=1}^K B \left( \vU_{h_{n_k}}(\vr_0, \vm_0) \right) \to \Ov{B(\vr, \vm)} \ \mbox{in}\
	L^1((0,T) \times \Td; R^{d+1}) \ \mbox{as}\ K \to \infty
\end{equation}	
a.s. for any Lipschitz $B: R^{d+1} \to R$.
The convergences \eqref{C12}, \eqref{C12a} hold for any subsequence of $(h_{n_k})_{k=1}^\infty$.
\end{Proposition}

To apply the above result to the Monte Carlo method, it is enough to realize that the mapping
\[
(\vr_0, \vm_0) \mapsto \frac{1}{K} \sum_{k=1}^K \vU_{h_{n_k}} (\vr_0, \vm_0)
\ \mbox{is Borel measurable}.
\]
Consequently, if $(\vr_{0,n}, \vm_{0,n})$ are pairwise i.i.d. data, then so are the corresponding approximations
\[
\frac{1}{K} \sum_{k=1}^K \vU_{h_{n_k}}(\vr_0, \vm_0) \ \mbox{as well as their limit}\ (\vr, \vm).
\]

\begin{Theorem}\label{TC1s} {\bf (Convergence of Monte Carlo method, strong form).}
	
	\noindent	
	Suppose $(\vr_{0,n}, \vm_{0,n})_{n=1}^\infty$ are a pairwise i.i.d.
	representations of random data $(\vr_0, \vm_0) \in \mathcal{D}$ satisfying \eqref{C1a} with a deterministic
	constant $L$. Let $(\vU_h)_{h \searrow 0}$ be a consistent approximation of the Euler system in the sense of Definition~\ref{D_CA}.
	
	Then any sequence $h_m \to 0$ contains a subsequence $(h_{m_k})_{k=1}^\infty$
	such that
	\begin{equation} \label{C13}
		\expe{ \left\| \frac{1}{N K} \sum_{n=1}^N
		\sum_{k=1}^K \vU_{h_{m_k}} (\vr_{0,n}, \vm_{0,n}) - \expe{ (\vr, \vm)   }\right \|_{L^q((0,T) \times \Td; R^{d+1})} }
		\to 0 \ \mbox{as}\ K,\ N \to \infty
	\end{equation}
for any $1 < q \leq \frac{2 \gamma}{\gamma + 1}$, 	
where $(\vr, \vm) \in \mathcal{U}(\vr_0, \vm_0)$ is a measurable
selection.
\end{Theorem}

Theorem \ref{TC1s} claims convergence of empirical averages in the strong Lebesgue topology $L^q$ to a measurable selection removing
both shortcomings of Theorem \ref{TC1}. The result holds for a suitable subsequence for which the  approximate
solutions (S)--converge  in the sense of \cite{Fei2021SC}.

\section{Convergence to smooth solutions}
\label{cs}

Finally, we rephrase Theorem \ref{TR1a} in terms of the Monte Carlo approximations. First observe that
any consistent approximation converges strongly in $L^q$ as long as the limit Euler system admits a
strong solution, see \cite[Chapter 7, Theorem 7.9]{FeLMMiSh}. Consequently, we obtain the following result.

	\begin{Theorem}\label{Tcs1} {\bf (Convergence of Monte Carlo method, smooth solutions).}
		
	\noindent	
	Suppose $(\vr_{0,n}, \vu_{0,n})_{n=1}^\infty$ is a pairwise independent
	representation of initial data $(\vr_0, \vm_0) \in \mathcal{D}_S$ satisfying
	\eqref{C1a} with a deterministic constant $L$. Let $(\vU_h)_{h \searrow 0}$ be a consistent approximation
	of the Euler system in the sense of Definition~\ref{D_CA}.

	Then
	\begin{align}
	\expe{ \left\| \frac{1}{N} \sum_{n=1}^N \mathds{1}_{t \leq T_{\rm max}} \vU_h (\vr_{0,n}, \vm_{0,n}) - \expe{ \mathds{1}_{t \leq T_{\rm max}} (\vr, \vm )}
	 \right\|_{L^q((0,T) \times \Td; R^{d + 1})} } \to 0, \quad 1 < q \leq \frac{ 2 \gamma}{\gamma + 1}
	 \nonumber
		\end{align}
as $h \to 0$, $N \to \infty$, where $(\vr, \vm)$ is the unique smooth solution of the
Euler system emanating from the data $(\vr_0, \vm_0)$.
	\end{Theorem}

\section{Monte Carlo finite volume method}
\label{FV}

The aim of this section is to illustrate theoretical results on a series of numerical simulations.
To this end, we firstly provide an example of a deterministic fully discrete numerical scheme for the Euler equations which yields consistent approximations according to Definition~\ref{D_CA}.

\subsection{Notations}

In what follows we briefly list all notations necessary to define our numerical scheme  based on the implicit time discretization and  finite volume method.

\subsubsection*{Space discretization}
Let $\grid_h$ be a structured mesh approximating the physical domain $\Td,$
\begin{align*}
\Td=\bigcup_{K\in\grid_h} K,
\end{align*}
where the elements $K$ are rectangles or cuboids. Let $h \in (0,1)$ be the  mesh parameter, meaning $|K| \approx h^d.$ Further, let $\mathcal{E}$ be the set of all faces $\sigma$ of elements $K,$ $|\sigma|\approx h^{d-1}.$
 By $Q_h$ we denote the space of piecewise constant functions on elements $K$,
\begin{equation*}
Q_h = \left\{\phi_h \in L^1(\Omega) \,|\, \phi_{h_{|_{K}}}=\phi_K = const. \mbox{ for all } K\in\grid_h \right\},
\end{equation*}
and by $\vc{Q}_h$ its analogue for vector-valued functions.
%There is a standard projection operator associated to this space given by
%\begin{equation*}
%\Pi_Q : L^1 (\Omega) \to Q_h, \quad \Pi_Q \phi =  %\sum_{K \in \grid} \frac{1_{K}}{|K|} \int_K \phi \dx.
%\end{equation*}
%Let $\phi_h \in Q_h$. Then $\phi_h^{\rm out},$ $\phi_h^{\rm in}$ stand for its outward and inward traces on a face $\sigma \in \mathcal{E}.$  We define the average and jump operators on any $\sigma$ by
%\begin{align*}
%\avs{\phi_h}=\frac{\phi_h^{\rm out}+\phi_h^{\rm in}}{2}, \qquad \jump{\phi_h}=\phi_h^{\rm out}-\phi_h^{\rm in},
%\end{align*}
%respectively.
For $\phi_h \in Q_h$ we denote the standard average and jump operators on any $\sigma$ by $\avs{\phi_h}$ and $\jump{\phi_h}$, respectively.
Further, we define
the discrete gradient and discrete divergence operators,
\begin{equation*}
\begin{aligned}
\Gradedge \phi_h  &=  \sum_{\sigma \in \mathcal{E}} \left( \Gradedge \phi_h\right)_\sigma  1_{\sigma}, &
\left( \Gradedge \phi_h\right)_{\sigma} & =  \frac{\jump{\phi_h}}{h}  \vn ,
\\
\Divh \bfpphi_h &= \sum_{K \in \grid_h} \left( \Divh \bfpphi_h \right)_K 1_K, &
\left( \Divh \bfpphi_h \right)_K &=  \sum_{\sigma \in \facesK} \frac{|\sigma|}{|K|}  \avs{\bfpphi_h} \cdot  \vn  ,
\end{aligned}
\end{equation*}
where  $\bfpphi_h=(\phi_{1,h}, \dots, \phi_{d,h})^T \in \vc{Q}_h,$ $\vc{n}$ is an unit outward normal vector to face $\sigma,$ $\mathcal{E}(K)$ is the set of all faces of an element $K,$ and $1_K,$ $1_{\sigma}$ are characteristic functions.

We consider the diffusive upwind numerical flux function of the form
\begin{align}
\Fup[r_h,\vc{u}_h]=\avs{r_h} \ \avs{\vc{u}_h} \cdot \vc{n} - \left(h^\varepsilon + \frac{1}{2} |\avs{\vc{u}_h} \cdot \vc{n}| \right)\jump{ r_h },  \qquad \varepsilon>-1.
\end{align}
$\vFup$ shall stand for its vector-valued analogue.   See \cite[Chapter 8]{FeLMMiSh} for more details.

\subsubsection*{Time discretization}
We consider an equidistant discretization of the time interval $[0,T]$ with the time step $\TS >0$ and time instances $t^k=k\TS,$ $k=0,\ldots,N_T.$
We denote $f^k(x) = f(t^k,x) $  for all $x\in \Td,$ by $f_\Delta$ a piecewise constant interpolation  of discrete values $f^k,$ and by $f_h^k$, $k=1,2,\ldots,N_T,$ a fully discrete function.  For its piecewise constant interpolation in time we use simplified notation $f_h$ instead of $f_{\Delta,h}$. The operator $D_t$ in \eqref{scheme} below stands for
\begin{align*}
 D_t f_h \equiv \frac{f_h(t) - f_h(t-\Delta t)}{\TS},
\end{align*}
which means we use the  backward Euler finite difference method to approximate the time derivative.

\subsection{Viscous finite volume (VFV) method}

Let the initial  data $(\vr_0,\vm_0) \in \mathcal{D}$ be given,  $(\vr_h^0,\vm_h^0) \in Q_h \times \vc{Q}_h$  be piecewise constant projections of the initial data, and $\vu_h^0 = \frac{\vm_h^0}{\vrh^0}.$
%\begin{align*}
%(\vr_h^0,\vm_h^0) = (\Pi_Q \vr_0, \Pi_Q \vm_0) \in Q_h \times \vc{Q}_h, \ \vu_h^0 = \frac{\vm_h^0}{\vrh^0}.
%\end{align*}
We seek piecewise constant functions (in space and  time) $(\vrh,\vuh) $   satisfying the following equations:
\begin{subequations}\label{scheme}
\begin{align}
\intTd{ D_t \vrh \varphi_h } &-   \sumG \intG{ \Fup [\vrh,\vuh]
\jump{\varphi_h}   } = 0 \quad \mbox{for all } \varphi_h \in Q_h,  \label{scheme1}
\\
\intTd{ D_t  (\vrh \vuh) \cdot \bfphi_h } &-   \sumG \intG{ \vFup [\vrh \vuh,\vuh]
\cdot \jump{\bfphi_h}   } - \intTd{  a\vrh^\gamma   \Divh \bfphi_h   } \nonumber
\\= &- h^\alpha \intTd{ \Gradedge \vuh   :  \Gradedge \bfphi_h  }
\quad \mbox{for all }
\bfphi_h \in \vc{Q}_h. \label{scheme2}
\end{align}
\end{subequations}
This finite volume scheme was introduced in \cite[Chapter 12]{FeLMMiSh}.

\begin{Remark}
Scheme \eqref{scheme} can be seen as a vanishing viscosity approximation of the barotropic Euler system due to a Newtonian-type viscosity of order $h^\alpha,$ $\alpha > 0$ included in the  momentum equation \eqref{scheme2}. This numerical diffusion term allows us to control the discrete velocity gradients.  As a result, the density remains strictly positive at any time level without imposing any extra CFL-type condition. Moreover, the VFV method is structure-preserving, meaning that it is conservative and dissipates discrete energy. These properties allow to show that the VFV method yields consistent approximations and consequently, its unconditional convergence, see \cite[Chapter 12]{FeLMMiSh}.
Note that standard finite volume schemes  usually require assumptions on uniform boundedness of numerical solutions in order to rigorously prove their consistency, see, e.g., \cite[Chapter 10]{FeLMMiSh}.
\end{Remark}

%
%Stability and consistency of the scheme stem from the following crucial properties:
%\begin{itemize}
%\item {\bf Existence of solutions:} There exist solutions $\{(\vrh^k, \vuh^k)\}_{k=1}^{N_T}$ to scheme \eqref{scheme}.
%\item {\bf Discrete conservation of mass:}
%$ \displaystyle
%	\intTd{ \vrh(t) } = \intTd{ \vr_0 }, \; t\in(0,T]
%$
%\item {\bf Positivity of discrete density:} $\vrh(t) > 0$ for all $t\in(0,T]$
%\item {\bf Discrete energy balance:}
%$ \displaystyle D_t \left(\intTd{ E(\vrh^k,\vmh^k)}\right)= e^k_{h},$ $e^k_{h} \leq 0$ for $k=1,\ldots, N_T$
%\end{itemize}
%We refer to \cite[Chapters 11-12, Lemmas 11.3 and 12.1]{FeLMMiSh} for more details and proofs.

\begin{Proposition}[\bf Consistency of the VFV method]\label{T_CA}
Suppose the initial data $(\vr_0, \vm_0) \in \mathcal{D}$ satisfy  \eqref{C1a} for some constant $L > 0$.
Let the parameters $\varepsilon$ and $\alpha$ satisfy
\begin{equation}\label{eps_alpha}
\begin{split}
& -1 < \varepsilon \
\mbox{ and } \ 0 < \alpha <  2 - \frac{ d/3 + 1 + \varepsilon }{\gamma} \quad \mbox{ for } \gamma\in(1,2),
\\
\mbox{or}\quad  & -1 < \varepsilon \ \mbox{ and } \ 0 < \alpha <  2 - \frac{ d}{\gamma} \quad \mbox{ for } \gamma \geq 2.  \\
\end{split}
\end{equation}

Then the numerical scheme \eqref{scheme} is consistent with the Euler equations \eqref{E1}-\eqref{E4} in the sense of Definition~\ref{D_CA}.
%  Specifically,  numerical solution $\vc{U}_h=\vc{U}_h(\vr_0,\vm_0)=(\vr_h, \vr_h \vu_h \equiv \vm_h)$ satisfies
%\begin{equation*}
%\int_0^T \intTd{ \left[ \vrh \partial_t \varphi + \vmh \cdot \Grad \varphi \right]} \dt  = \intTd{ \vrh^0 \varphi(0,\cdot) }  + \int_0^T
%e_{1,h} (t, \varphi) \dt
%\end{equation*}
%for any $\varphi \in C_c^2([0,T) \times \Td)$;
%\begin{eqnarray*}
% &&
%\int_0^T \intTd{ \left[ \vmh \cdot \partial_t \boldsymbol{\varphi} + \left(\frac{\vmh \otimes \vmh}{\vrh} \right) : \Grad \boldsymbol{\varphi}  + a\vrh^\gamma \Div \boldsymbol{\varphi} \right]} \dt
%\nonumber \\&&
% =  \intTd{ \vmh^0 \cdot \boldsymbol{\varphi}(0,\cdot) }   + \int_0^T e_{2,h} (t, \boldsymbol{\varphi}) \dt,
%\end{eqnarray*}
%for any $\boldsymbol{\varphi} \in C^2_c([0,T) \times \Td; R^d)$;
%\begin{equation*}
%	\intTd{ E(\vr_h, \vm_h) (\tau, \cdot) } \leq \intTd{ E(\vrh^0, \vmh^0) } +  e_{3,h}
%	\end{equation*}
%for a.a. $\tau \in (0,T)$;
%\[
%\| e_{1,h} (\cdot, \varphi ) \|_{L^1(0,T)} \leq c(L) h^{\beta}  \| \varphi \|_{C^2}, \ \
%\| e_{2,h} (\cdot, \boldsymbol{\varphi} ) \|_{L^1(0,T)} \leq c(L) h^{\beta}  \| \boldsymbol{\varphi} \|_{C^2}\, \mbox{ for some }\ \beta > 0,
%\]
%and $\displaystyle e_{3,h}=\int_0^\tau e_h(t) \, \dt  \leq 0.$
\end{Proposition}
%
%Note that $e^1(h,L,\|\varphi\|_{C^2})=\| e_{1,h} (\cdot, \varphi ) \|_{L^1(0,T)}$ and  $e^2(h,L,\|\bfphi\|_{C^2})=\| e_{2,h} (\cdot, \bfphi ) \|_{L^1(0,T)}$  obviously vanish as $h\rightarrow 0.$

%Assuming the Euler system admits a strong solution, finite volume approximations generated by scheme \eqref{scheme} are unconditionally convergent in the sense that no assumption on boundedness of numerical solutions is required.
Since the finite volume scheme \eqref{scheme} generates consistent approximations of the Euler system, the main results on the convergence of  the Monte Carlo method are, without  further assumptions, valid for numerical solutions of the Euler system computed by the scheme \eqref{scheme}.

\begin{Corollary}
Theorems~\ref{TC1}, \ref{TC1s} and \ref{Tcs1} remain valid for numerical approximations $(\vc{U}_h)_{h \searrow 0}$  obtained as solutions to the finite volume scheme \eqref{scheme} approximating the Euler system \eqref{E1}-\eqref{E4} provided the parameters $\varepsilon$ and $\alpha$ satisfy \eqref{eps_alpha}.
\end{Corollary}

\subsection{Monte Carlo VFV simulations}

We consider two-dimensional Kelvin-Helmholtz problem on a computational domain $[0,1]\times [0,1]$ with
periodic boundary conditions and
the initial data
\begin{equation}
	\left(\varrho, u_1, u_2\right)(x,0)= \begin{cases}(2,-0.5,0), & \text { if } I_1<x_2<I_2 \\ (1,0.5,0), & \text { otherwise. }\end{cases}
\end{equation}
The interface profiles  $I_1, I_2$ are random variables given as
\begin{equation}
	I_j(x, \omega):=J_j+\epsilon Y_j(x, \omega), \quad j=1,2, \,\epsilon = 0.01,
\end{equation}
and
\begin{equation}
	J_1=0.25, J_2=0.75, Y_j(x, \omega)=\sum_{i=1}^m a_j^i(\omega) \cos \left(b_j^i(\omega)+2 i \pi x_1\right), \quad j=1,2, \ m=10.
\end{equation}
Here $a_j^i(\omega) \sim {\Bbb U }([0,1])$ and $b_j^i(\omega)\sim {\Bbb U}([-\pi,\pi]),$ $j=1,2,$ are uniformly distributed random variables.  The coefficients $a_{j}^i$ have been normalized such that $\sum_{i=1}^{m}a_j^{i}=1, \, j=1,2.$

In the VFV method \eqref{scheme} the parameters of numerical flux are chosen to be $\alpha=0.8, \varepsilon=0$ to satisfy condition \eqref{eps_alpha} (note that $\gamma =1.4.$)  Let $(\varrho_{0,n}, \bm{u}_{0,n})$, $n=1,2, \dots, N,$ be pairwise
i.i.d. samples from the random initial data $(\varrho_{0}, \bm{u}_{0}).$  The   VFV solutions corresponding to the initial data $(\varrho_{0,n}, \bm{u}_{0,n})$
obtained on a mesh $\mathcal{T}_h$
are denoted by $\bm{U}^n_h,$ $n=1, \dots, N.$

\begin{figure}[!h]
	\centering
	\subfigure{
		\includegraphics[width=0.4\linewidth]{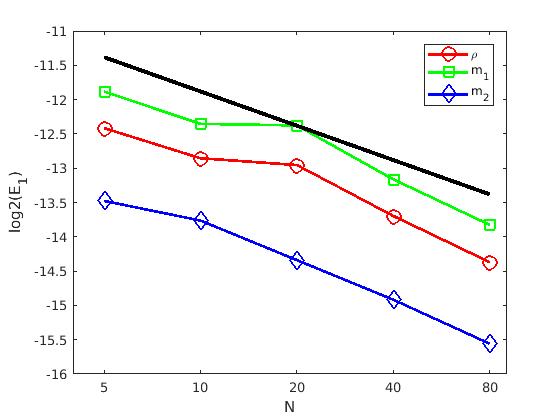}
	}
	\subfigure{
		\includegraphics[width=0.4\linewidth]{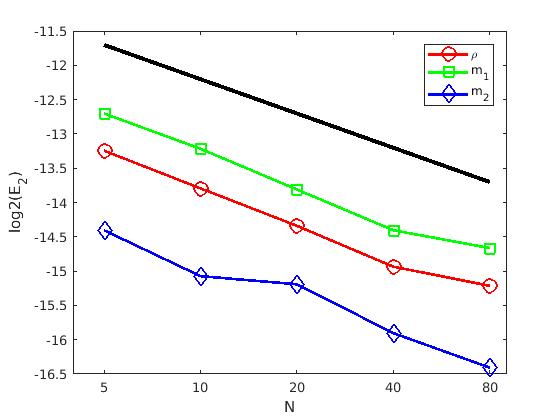}
	}
	\caption{The statistical errors  $E_1$ (left)  and  $E_2$(right) computed by the Monte Carlo VFV method on the mesh with $128 \times 128$ cells. The black solid lines without any marker denote the reference slope of $N^{-1/2}$.}
	\label{fig: convergence-128}
\end{figure}

\begin{figure}[!h]
	\centering
	\subfigure{
		\includegraphics[width=0.4\linewidth]{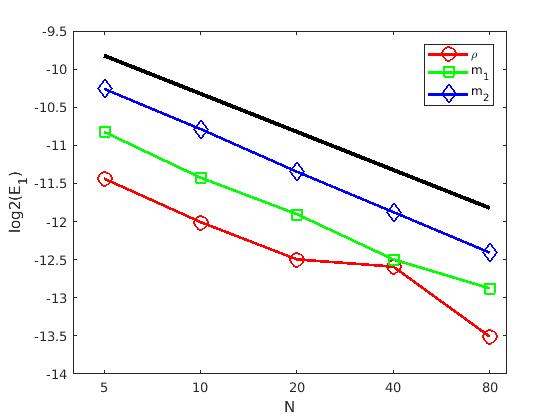}
	}
	\subfigure{
		\includegraphics[width=0.4\linewidth]{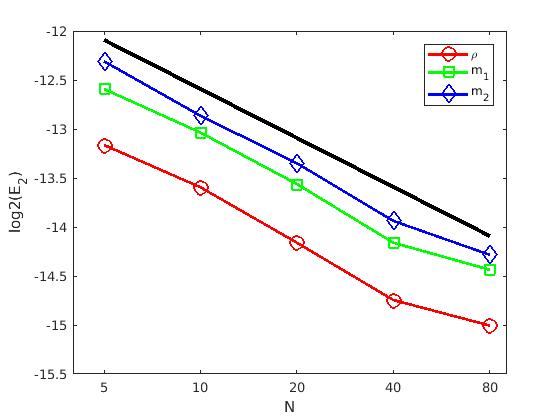}
	}
	\caption{The statistical errors  $E_1$ (left)  and  $E_2$(right) computed by the Monte Carlo VFV method on the mesh with $512 \times 512$ cells. The black solid lines without any marker denote the reference slope of $N^{-1/2}$.}
	\label{fig: convergence-512}
\end{figure}

\begin{figure}[!h]
	\centering
	\subfigure{
		\includegraphics[width=0.4\linewidth]{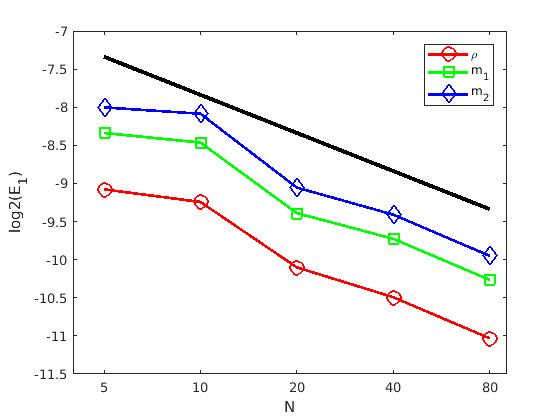}
	}
	\subfigure{
		\includegraphics[width=0.4\linewidth]{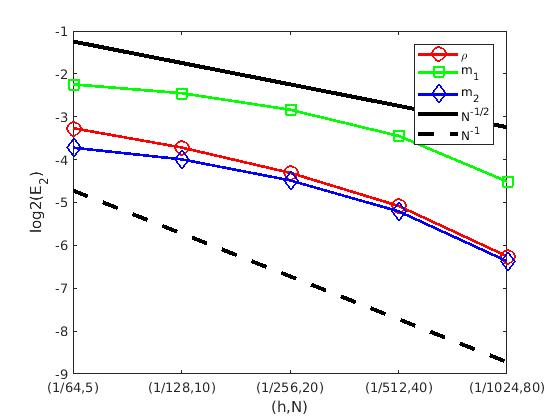}
	}
	\caption{The statistical errors  $E_1$ (left)  and  $E_2$(right) computed by the Monte Carlo VFV method on the mesh with $2048 \times 2048$ cells. The black solid lines without any marker denote the reference slope of $N^{-1/2}$.}
	\label{fig: convergence-2048}
\end{figure}

Our aim is to demonstrate the convergence of the Monte Carlo VFV method in its weak and strong form, cf.~Theorem~\ref{TC1} and Theorem~\ref{TC1s}. The latter considers $\mathcal{K}$-convergence of VFV solutions, i.e.~convergence of the Ces\'aro averages with respect to different meshes $\mathcal{T}_{h_{n_k}}$ that are computed as
$$
\widetilde{\bm{U}}_{h_{K}} := \frac{1}{K} \sum_{k=1}^K \bm{U}_{h_{n_k}}.
$$
Clearly, a projection operator between coarse and fine grids has to be applied to compute the above sum over different mesh resolutions.
We consider the following error functions
\begin{equation}
\label{eq7}
E_1(\bm{U}_h^N)=\frac 1 L \sum_{\ell =1}^L\left\|\frac{1}{N} \sum_{n=1}^N\left(\bm{U}^{n, \ell}_{h}(T, \cdot)-\mathbb{E}[\bm{U}_{h}(T, \cdot)]\right)\right\|_{L^1\left(\mathbb{T}^d\right)},
\end{equation}
\begin{equation}
\label{eq7a}
E_2(\bm{U}^N_{h_K})=\frac 1 L \sum_{\ell =1}^L\left\|\frac{1}{N} \sum_{n=1}^N\left(\widetilde{\bm{U}}_{h_{K}}^{n, \ell}(T, \cdot)-\mathbb{E}[\widetilde{\bm{U}}_{h_{K}}(T, \cdot)]\right)\right\|_{L^1\left(\mathbb{T}^d\right)}
\end{equation}
with $N=5,10,20,40,80$ and $L=20.$
We work with consecutively refined meshes with mesh parameters $h_1=1 /64$, $h_2=1/128, \dots, h_6 = 1/2048.$
%The reference grid parameters are $h_{ref} = 1/2048$ and $K_{ref} = 6.$
The reference expected values $\mathbb{E}[\bm{U}_{h}], \, \mathbb{E}[\widetilde{\bm{U}}_{h_{K}}] $ are computed using $N_{ref}=100$ samples
\begin{equation}
\label{eq9}
	\mathbb{E}[\bm{U}_{h}]=\frac{1}{N_{ref}} \sum_{s=1}^{N_{ref}} \bm{U}_{h}^s,\quad \mathbb{E}[\widetilde{\bm{U}}_{h_{K}}]
=\frac{1}{N_{ref}} \sum_{s=1}^{N_{ref}} \widetilde{\bm{U}}_{h_{K}}^s,\quad N_{ref}=100.
\end{equation}
Behaviour of statistical errors $E_1$ and $E_2$ on a  reference mesh with $h_{ref}=2048$ and $K_{ref}=6$ are shown in Tables~\ref{E1: Barotropic-2048}, Table~\ref{E2: Barotropic-2048}. Note that formulae \eqref{eq7}--\eqref{eq9} are applied componentwise for $\vr, m_1 = \varrho u_1,  m_2 = \varrho u_2.$ Figures~\ref{fig: convergence-128}, \ref{fig: convergence-512} and
\ref{fig: convergence-2048} illustrate the convergence with respect to $N$ random samples on meshes with $128\times 128,$ $512 \times 512$ and $2048\times2048$ cells, respectively. We can observe that applying  the C\'esaro averages statistical convergence rate of $N^{-1/2}$ is obtained already on coarser meshes. Although in Theorem~\ref{TC1} the  convergence is rigorously proved  only in the weak topology our numerical results for $E_1$ indicate the convergence even in  $L^1$-norm.

\begin{figure}[!h]
	\centering
	\subfigure[$E_1$]{
		\includegraphics[width=0.4\linewidth]{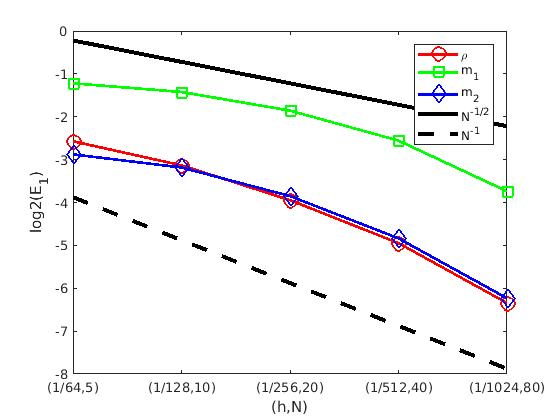}
	}
	\subfigure[$E_2$]{
		\includegraphics[width=0.4\linewidth]{pictures/convergence-couple-E1-2.jpg}
	}
	\caption{The errors $E_1$ and $E_2$ in $L^1-$norm with parmaters $h=1/(2^{\ell+5}), N(h)=5\cdot 2^{\ell-1}, \ell=1,2,\dots,5.$ The black solid lines without any marker denote the reference slope of $N^{-1/2}$. The black dash lines without any marker denote the reference slope of $N^{-1}$.}
	\label{FIG_total_error}
\end{figure}

Finally, in order to consider both the approximation and statistical errors, we present in Table~\ref{E1: Barotropic_total_error}, Table~\ref{E2: Barotropic_total_error} and
Figure~\ref{FIG_total_error} the
total errors $E_i(U_h^N(h)), i=1,2,$ with respect to the parameters pair $(h, N(h))$, $h=1/ (2^{\ell +5});
N(h) = 5 \cdot 2^{\ell -1}, \quad \ell =1,\dots,5.$ In formulae \eqref{eq7}, \eqref{eq7a} the reference expected values were computed using $h_{ref}=1/2048$ and $N_{ref} = 100.$ As illustrated in
Table~\ref{E1: Barotropic_total_error}, Table~\ref{E2: Barotropic_total_error} and
Figure~\ref{FIG_total_error}
the convergence rates are between $N^{-1/2}$ and $N^{-1}$.

\begin{table}[!h]
	\centering
	\caption{The statistical errors $E_1$ computed by the Monte Carlo VFV method on  the mesh with $2048\times 2048$ cells.}
\medskip
	\begin{tabular}{|c|cc|cc|cc|}
		\hline
		{variables} & \multicolumn{2}{c|}{$\varrho$} & \multicolumn{2}{c|}{$m_1$}  & \multicolumn{2}{c|}{$m_2$}
\\ \hline
		samples& error  & order & error  & order  & error  & order     \\
		\hline
		5  & 1.85e-03      &          & 3.09e-03   &      &3.90e-03& \\
		10  &1.65e-03    & 0.17  &2.83e-03   & 0.13  &3.68e-03&0.08 \\
		20  &9.10e-04    & 0.86 & 1.49e-03  & 0.93  &1.88e-03 &0.97\\
		40  &6.94e-04  & 0.39  & 1.18e-03   & 0.34 &  1.47e-03&0.35\\
		80 &4.77e-04  & 0.54 &8.12e-04  & 0.54  &  1.01e-03&0.54\\
		\hline
	\end{tabular}
	\label{E1: Barotropic-2048}
\end{table}
\begin{table}[!h]
	\centering
	\caption{The statistical errors $E_2$ computed by the Monte Carlo VFV method on  the mesh  with $2048\times 2048$ cells.}
\medskip
	\begin{tabular}{|c|cc|cc|cc|}
		\hline
		{variables} & \multicolumn{2}{c|}{$\varrho$} & \multicolumn{2}{c|}{$m_1$}  & \multicolumn{2}{c|}{$m_2$}   \\ \hline
		samples& error  & order & error  & order  & error  & order     \\
		\hline
		5  & 3.97e-04      &         & 6.64e-04   &      &8.76e-04& \\
		10  &2.81e-04    & 0.50  &4.65e-04   & 0.51  &6.03e-04&0.54 \\
		20  &1.89e-04    & 0.57 & 3.14e-04  & 0.57  &3.95e-04 &0.61\\
		40  &1.32e-04  & 0.52  & 2.21e-04   & 0.51 &  2.87e-04&0.46\\
		80 &1.08e-05  & 0.29 &1.79e-04  & 0.30 &  2.27e-04&0.34\\
		\hline
	\end{tabular}
	\label{E2: Barotropic-2048}
\end{table}

\begin{table}[!h]
	\centering
	\caption{The total error $E_1$ computed by the Monte Carlo VFV method with parameters $(h,N(h))$; \ $h=1/ (2^{\ell +5}), N(h) = 5 \cdot 2^{\ell -1}, \ \ell =1,\dots,5.$ }
\medskip
	\begin{tabular}{|c|cc|cc|cc|}
		\hline
		{variables} & \multicolumn{2}{c|}{$\varrho$} & \multicolumn{2}{c|}{$m_1$}  & \multicolumn{2}{c|}{$m_2$} \\ \hline
		$(h,N(h))$ & error  & order & error  & order  & error  & order      \\
		\hline
		$(1/64,5)$ & 1.68e-01     &         & 4.29e-01   &     &1.36e-01& \\
		$(1/128,10)$ &1.14e-01    & 0.56  &3.73e-01   & 0.20  &1.10e-01&0.31 \\
		$(1/256,20)$ &6.49e-02    & 0.81 & 2.76e-01  & 0.43  &6.93e-02 &0.67\\
		$(1/512,40)$ &3.24e-02  & 1.00  & 1.70e-01   & 0.70 &  3.50e-02&0.99\\
		$(1/1024,80)$ &1.23e-02  & 1.40 &7.49e-02  & 1.18&  1.33e-02&1.40\\
		\hline
	\end{tabular}
	\label{E1: Barotropic_total_error}
\end{table}

\begin{table}[!h]
	\centering
	\caption{The total error $E_2$ computed by the Monte Carlo VFV method with  parameters $(h,N(h))$; \ $h=1/ (2^{\ell +5}), N(h) = 5 \cdot 2^{\ell -1}, \ \ell =1,\dots,5.$ }
\medskip
	\begin{tabular}{|c|cc|cc|cc|}
		\hline
		{variables} & \multicolumn{2}{c|}{$\varrho$} & \multicolumn{2}{c|}{$m_1$}  & \multicolumn{2}{c|}{$m_2$} \\ \hline
		$(h,N(h))$ & error  & order & error  & order  & error  & order      \\
		\hline
		$(1/64,5)$ & 1.04e-01     &         & 2.11e-01   &     &7.58e-02& \\
		$(1/128,10)$ &7.60e-02   & 0.45  &1.83e-01   & 0.21  &6.28e-02&0.27 \\
		$(1/256,20)$ &5.05e-02    & 0.59 & 1.40e-01  & 0.39  &4.47e-02 &0.49\\
		$(1/512,40)$ &2.96e-02  & 0.77  & 9.17e-02   & 0.61 &  2.71e-02&0.72 \\
		$(1/1024,80)$ &1.31e-02  & 1.18 &4.39e-02  & 1.06 &  1.21e-02&1.16\\
		\hline
	\end{tabular}
	\label{E2: Barotropic_total_error}
\end{table}

\newpage

\section*{Conclusions}

In the present paper we have studied convergence of the Monte Carlo method combined with a consistent approximation for the random isentropic Euler equations. We work with the concept of dissipative weak solutions that
can be seen as a universal closure of consistent approximations.  Since the Euler equations are not uniquely solvable in the class of dissipative weak solutions we apply
the set-valued version of the Strong law of large numbers for general multivalued mapping with
closed range, cf.~Theorem~\ref{TR1}. For strong solutions this reduces  to the strong convergence in $L^q((0,T)\times T^d),$ $q= \frac{2 \gamma}{\gamma +1}.$  Combining Theorem~\ref{TR1} with the deterministic convergence results of a consistent approximation yield the convergence of the Monte Carlo method in the weak form, cf.~Theorem~\ref{TC1}. Applying
$\mathcal{K}-$convergence in random setting we have derived the  convergence of the Monte Carlo method in the strong form in Theorem~\ref{TC1s}. If the strong solution exists,  we obtain the strong convergence of the Monte Carlo estimators obtained by the numerical approximation  to the expected value of the unique strong statistical solution, cf.~Theorem~\ref{Tcs1}. In Section~\ref{FV} we illustrate theoretical results by numerical simulations obtained by the
 Monte Carlo method combined with the viscosity finite volume method.

%\bibliographystyle{plain}
%\bibliography{citace}

\def\cprime{$'$} \def\ocirc#1{\ifmmode\setbox0=\hbox{$#1$}\dimen0=\ht0
	\advance\dimen0 by1pt\rlap{\hbox to\wd0{\hss\raise\dimen0
			\hbox{\hskip.2em$\scriptscriptstyle\circ$}\hss}}#1\else {\accent"17 #1}\fi}

\end{document}